\documentclass{article}
\usepackage{amssymb}
\usepackage{amsmath}
\usepackage{amsfonts}
\usepackage{amsthm}
\usepackage{graphicx}

\newtheorem{theorem}{Theorem}

\newtheorem{definition}[theorem]{Definition}
\newtheorem{example}[theorem]{Example}

\newtheorem{lemma}[theorem]{Lemma}

\newtheorem{proposition}[theorem]{Proposition}
\newtheorem{remark}[theorem]{Remark}

\newcommand{\baseRing}[1]{\ensuremath{\mathbb{#1}}}

\newcommand{\N}{\baseRing{N}}
\newcommand{\R}{\baseRing{R}}
\newcommand{\Rs}{\overline\R}

\newcommand{\e}{\varepsilon}

\newcommand{\FIXME}[1]{???}

\DeclareMathOperator{\Id}{Id}

\newcommand{\dsum}{\sum}
\newcommand{\diint}{\iint}
\newcommand{\dint}{\int}
\newcommand{\tint}{\int}

\begin{document}

\title{Solving a characteristic Cauchy problem }
\author{Emmanuel ALLAUD, Victor DEVOUE \\
Equipe Analyse Alg\'{e}brique Non Lin\'{e}aire\\
D\'{e}partement Scientifique Interfacultaire\\
Universit\'{e} des Antilles et de la Guyane\\
Campus de Schoelcher. BP 7209\\
97275 Schoelcher Cedex, Martinique. F.W.I.}
\date{10 September, 2008}
\maketitle

\begin{abstract}
In this paper we give a meaning to the nonlinear characteristic Cauchy
problem for the Wave Equation in base form by replacing it by a family of
non-characteristic problems in an appropriate algebra of generalized
functions. We prove existence of a solution and we precise how it depends on the choice made. We also check that in the classical case (non-characteristic) our new solution coincides with the classical one.

\textit{MSC: 35D05; 35L05; 35L70; 46F30.}

\textbf{Key words : }algebras of generalized functions; nonlinear partial
differential equations; Characteristic Cauchy problem; Wave Equation.
\end{abstract}

\section{Introduction}

The goal of this paper is to study the solution to the following
characteristic Cauchy problem%
\begin{equation*}
\left( P_{C}\right) \left\{ 
\begin{array}{l}
\dfrac{\partial ^{2}u}{\partial x\partial y}=F(.,.,u), \\ 
u_{\mid \gamma }=\varphi , \\ 
\dfrac{\partial u}{\partial y}\mid _{\gamma }=\psi \text{,}%
\end{array}%
\right. 
\end{equation*}%
where $\gamma $ is a curve supposed to be a \textrm{C}$^{\infty }$ manifold
characteristic for the problem, $\varphi $ and $\psi $ being smooth
functions defined in $\gamma $.\textit{\ }The initial values have to be
defined as restrictions of functions (even in a generalized sense) given in $%
\gamma $\textit{\ }whose equation is $y=f(x)$.

As there is no classical solution here, we will look for a solution in a
broader context, using the framework of generalized functions \cite{Bia}, 
\cite{GKOS}, \cite{NePiScar}, \cite{NOP}. They are an efficient tool to
solve nonlinear problems as in \cite{JAM4}, \cite{MNV}. The general idea
goes as follows. The characteristic problem is approached by a one-parameter
family of classical smooth problems by deforming the characteristic curve $%
y=f(x)$ into a family of non-characteristic ones $y=f_{\varepsilon }(x)$; we
then get a one-parameter family of classical solutions. That is where the
framework of generalized functions is used ; by means of this
regularization, we define an associated generalized problem and we interpret
this family of solutions as a generalized solution itself. Indeed a
generalized function can be defined as a one-parameter family of smooth
functions satisfying some asymptotical growth restrictions \cite{NOP}.

More precisely we will take $(f_{\varepsilon })_{\varepsilon }$ to be
equivalent to $f$ for some sense in an appropriate algebra of generalized
functions. Furthermore, by imposing some restrictions on the asymptotical
growth of the $f_{\varepsilon }$, we are able to prove that the generalized
solution depends solely on the class of $(f_{\varepsilon })_{\varepsilon }$
as a generalized function, not on the particular representative. We also
prove that in the non-characteristic smooth case, the generalized solution
provided by our method coincides (in the sense of generalized functions)
with the classical smooth solution.

The plan of this article is as follows. This section is followed by section %
\ref{sect_defn} which briefly introduces the generalized algebras with our
application in mind. We define a generalized differential problem associated
to the ill posed classical one, then we proceed in section \ref%
{sect_existence} with the proof of the existence of the generalized
solution. In our case this amounts to prove that, provided a set of
restrictions on the curve and the chosen deformation, the one-parameter
family of solutions satisfy the required asymptotical growth. Subsection \ref%
{sect_unicite} is devoted to prove that the generalized solution does not
depend on the representative of $(\gamma _{\varepsilon })_{\varepsilon }$.

Then in section \ref{sect_exemples} we compute a few examples of
characteristic equations, and make the link with distributional solutions.

\section{Algebras of generalized functions}

\label{sect_defn} We recall briefly here the definition of the $(\mathcal{C},%
\mathcal{E},\mathcal{P})$-algebras upon which the remaining of this paper is
based. We follow here the expositions found in \cite{JAM1}, \cite{JAM2}, 
\cite{JAM3}, \cite{JAM4}, \cite{MNV} and \cite{Dev1}. We do not intend to
properly define and explain $(\mathcal{C},\mathcal{E},\mathcal{P})$-algebras
here, we rather want to fix the notations we will use in the latter
sections. We refer the reader to the references.

The formalism described here, the $(\mathcal{C},\mathcal{E},\mathcal{P})$%
-algebras, is well suited for partial differential equations because of its
parametric nature; this will become clear in the next section.

\subsection{The presheaves of $(\mathcal{C},\mathcal{E},\mathcal{P})$%
-algebras}

\subsubsection{Definitions}

Take

\begin{itemize}
\item $\Lambda$ a set of indices;

\item $A$ a solid subring of the ring $\mathbb{K}^{\Lambda}$,\ ($\mathbb{K=R}
$ or $\mathbb{C}$), that is $A$ has the following stability property:
whenever $(\left\vert s_{\lambda}\right\vert
)_{\lambda}\leq(r_{\lambda})_{\lambda}$ (i.e. for any $\lambda$, $\left\vert
s_{\lambda}\right\vert \leq r_{\lambda}$) for any pair $\left(
(s_{\lambda})_{\lambda},(r_{\lambda})_{\lambda}\right) \in\mathbb{K}%
^{\Lambda}\times\left\vert A\right\vert $, it follows that $%
(s_{\lambda})_{\lambda}\in A$, with $\left\vert A\right\vert =\{(\left\vert
r_{\lambda}\right\vert )_{\lambda}:(r_{\lambda})_{\lambda}\in A\}$;

\item $I_{A}$ an solid ideal of $A$ with the same property;

\item $\mathcal{E}$ a sheaf of $\mathbb{K}$-topological algebras on a
topological space $X$, such that for any open set $\Omega$ in $X$, the
algebra $\mathcal{E}(\Omega)$ is endowed with a family $\mathcal{P}(\Omega
)=(p_{i})_{i\in I(\Omega)}$ of seminorms satisfying 
\begin{equation*}
\forall i\in I(\Omega)\text{, }\exists(j,k,C)\in I(\Omega)\times
I(\Omega)\times\mathbb{R}_{+}^{\ast}\text{, }\forall f,g\in\mathcal{E}%
(\Omega):p_{i}(fg)\leq Cp_{j}(f)p_{k}(g).
\end{equation*}
Assume that

\item For any two open subsets $\Omega_{1}$, $\Omega_{2}$ of $X$ such that $%
\Omega_{1}\subset$ $\Omega_{2},$ we have $I(\Omega_{1})\subset I(\Omega_{2})$
and if $\rho_{1}^{2}$ is the restriction operator $\mathcal{E}(\Omega
_{2})\rightarrow\mathcal{E}(\Omega_{1})$, then, for each $p_{i}\in \mathcal{P%
}(\Omega_{1})$, the seminorm $\widetilde{p}_{i}=p_{i}\circ\rho _{1}^{2}$
extends $p_{i}$ to $\mathcal{P}(\Omega_{2})$;

\item For any family $\mathcal{F}=(\Omega _{h})_{h\in H}$ of open subsets of 
$X$ if $\Omega =\cup _{h\in H}\Omega _{h}$, then, for each $p_{i}\in 
\mathcal{P}(\Omega )$, $i\in I(\Omega )$, there exists a finite subfamily $%
\Omega _{1},\cdots,\Omega _{n\left( i\right) }$ of $\mathcal{F}$ and
corresponding seminorms $p_{1}\in \mathcal{P}(\Omega _{1}),\cdots,p_{n\left(
i\right) }\in \mathcal{P}(\Omega _{n\left( i\right) })$, such that, for each 
$u\in \mathcal{E}(\Omega )$,%
\begin{equation*}
p_{i}\left( u\right) \leq p_{1}\left( u_{\left\vert \Omega _{1}\right.
}\right) +\cdots+p_{n\left( i\right) }(u_{\left\vert \Omega _{n\left( i\right)
}\right. }).
\end{equation*}%
Set%
\begin{align*}
\mathcal{X}_{(A,\mathcal{E},\mathcal{P})}(\Omega )& =\{(u_{\lambda
})_{\lambda }\in \left[ \mathcal{E}(\Omega )\right] ^{\Lambda }:\forall i\in
I(\Omega )\text{, }\left( (p_{i}(u_{\lambda })\right) _{\lambda }\in
\left\vert A\right\vert \}, \\
\mathcal{N}_{(I_{A},\mathcal{E},\mathcal{P})}(\Omega )& =\{(u_{\lambda
})_{\lambda }\in \left[ \mathcal{E}(\Omega )\right] ^{\Lambda }:\forall i\in
I(\Omega )\text{, }\left( p_{i}(u_{\lambda })\right) _{\lambda }\in
\left\vert I_{A}\right\vert \}, \\
\mathcal{C}& =A/I_{A}.
\end{align*}%
One can prove that $\mathcal{X}_{(A,\mathcal{E},\mathcal{P})}$\textit{\ }is
a sheaf of subalgebras of the sheaf $\mathcal{E}^{\Lambda }$\textit{\ }and%
\textit{\ }$\mathcal{N}_{(I_{A},\mathcal{E},\mathcal{P})}$ is a sheaf of
ideals of $\mathcal{X}_{(A,\mathcal{E},\mathcal{P})}$ \cite{JAM2}. Moreover,
the constant sheaf $\mathcal{X}_{(A,\mathbb{K},\left\vert .\right\vert )}/%
\mathcal{N}_{(I_{A},\mathbb{K},\left\vert .\right\vert )}$\ is exactly the
sheaf\textit{\ }$\mathcal{C}=A/I_{A}$, and if $\mathbb{K}=\R$, $\mathcal C$ will be denoted $\Rs$.
\end{itemize}

\begin{definition}
We call presheaf of $(\mathcal{C},\mathcal{E},\mathcal{P})$-algebra the
factor presheaf of algebras over the ring $\mathcal{C}=A/I_{A}$%
\begin{equation*}
\mathcal{A}=\mathcal{X}_{(A,\mathcal{E}{\LARGE ,}\mathcal{P})}/\mathcal{N}%
_{(I_{A},\mathcal{E}{\LARGE ,}\mathcal{P})}.
\end{equation*}%
We denote by $\left[ u_{\lambda }\right] $\ the class in $\mathcal{A}(\Omega
)$ defined by the representative\textit{\ }$(u_{\lambda })_{\lambda \in
\Lambda }\in \mathcal{X}_{(A,\mathcal{E},\mathcal{P})}(\Omega ).$
\end{definition}

\subsubsection{\textbf{Overgenerated rings}}

See \cite{JAMADEV2}. \textit{Let }$B_{p}=\left\{ \left( r_{n,\lambda
}\right) _{\lambda }\in (\mathbb{R}_{+}^{\ast })^{\Lambda
}:n=1,\cdots,p\right\} $ \textit{and} $B$\textit{\ be the subset of }$(\mathbb{R%
}_{+}^{\ast })^{\Lambda }$\textit{\ } obtained as rational functions with
coefficients in $\mathbb{R}_{+}^{\ast }$, of elements in $B_{p}$ as
variables. Define 
\begin{equation*}
A=\left\{ \left( a_{\lambda }\right) _{\lambda }\in \mathbb{K}^{\Lambda
}\mid \exists \left( b_{\lambda }\right) _{\lambda }\in B,\exists \lambda
_{0}\in \Lambda ,\forall \lambda \prec \lambda _{0}:\left\vert a_{\lambda
}\right\vert \leq b_{\lambda }\right\} .
\end{equation*}

\begin{definition}
\label{Overgenerated}In the above situation, we say that $A$ is \emph{%
overgenerated }by $B_{p}$ (and it is easy to see that $A$ is a solid subring
of $\mathbb{K}^{\Lambda }$). If $I_{A}$\ is some solid ideal of $A$, we also
say that\ $\mathcal{C}=A/I_{A}$ is \emph{overgenerated }by $B_{p}$.
\end{definition}

\begin{example}
For example, as a \textquotedblleft canonical\textquotedblright\ ideal of $A$%
, we can take 
\begin{equation*}
I_{A}=\left\{ \left( a_{\lambda}\right) _{\lambda}\in\mathbb{K}^{\Lambda
}\mid\forall\left( b_{\lambda}\right) _{\lambda}\in
B,\exists\lambda_{0}\in\Lambda,\forall\lambda\prec\lambda_{0}:\left\vert
a_{\lambda}\right\vert \leq b_{\lambda}\right\} \text{.}
\end{equation*}

\begin{remark}
We can see that with this definition $B$ is stable by inverse.
\end{remark}
\end{example}

\subsubsection{Relationship with distribution theory}

Let $\Omega $ an open subset of $\mathbb{R}^{n}$. The space of distributions 
$\mathcal{D}^{\prime }(\Omega )$ can be embedded into $\mathcal{A}(\Omega )$. If $\left( \varphi _{\lambda }\right) _{\lambda \in (0,1]}$ is a family
of mollifiers $\varphi _{\lambda }\left( x\right) =\frac{1}{\lambda ^{n}}%
\varphi \left( \frac{x}{\lambda }\right) $, $x\in \mathbb{R}^{n}$, $\int
\varphi \left( x\right) dx=1$ and if $T\in \mathcal{D}^{\prime }\left( 
\mathbb{R}^{n}\right) $, the convolution product family $\left( T\ast
\varphi _{\lambda }\right) _{\lambda }$ is a family of smooth functions
slowly increasing in $\frac{1}{\lambda }$. So, for $\Lambda =(0,1]$, we
shall choose the subring $A$ overgenerated by some $B_{p}$ of $(\mathbb{R}%
_{+}^{\ast })^{\Lambda }$ containing the family $\left( \lambda \right)
_{\lambda }$, \cite{Del1}, \cite{Ober1}.

\subsubsection{\textbf{The association process}}

We assume that $\Lambda $ is left-filtering for a given partial order
relation $\prec $. We denote by $\Omega $ an open subset of $X$, $E$ a given
sheaf of topological $\mathbb{K}$-vector spaces containing $\mathcal{E}$ as
a subsheaf, $a$ a given map from $\Lambda $ to $\mathbb{K}$ such that $%
(a\left( \lambda \right) )_{\lambda }=(a_{\lambda })_{\lambda }$ is an
element of $A$. We also assume that%
\begin{equation*}
\mathcal{N}_{(I_{A},\mathcal{E}{\LARGE ,}\mathcal{P})}(\Omega )\subset
\left\{ \left( u_{\lambda }\right) _{\lambda }\in \mathcal{X}_{(A,\mathcal{E}%
{\LARGE ,}\mathcal{P})}(\Omega ):\lim\limits_{E(\Omega ),\Lambda }u_{\lambda
}=0\right\} .
\end{equation*}

\begin{definition}
We say that $u=\left[ u_{\lambda}\right] $ and $v=\left[ v_{\lambda }\right]
\in\mathcal{E}(\Omega)$ are $a$-$E$ associated if%
\begin{equation*}
\lim\limits_{E(\Omega),\Lambda}a_{\lambda}(u_{\lambda}-v_{\lambda})=0\text{.}
\end{equation*}
\textit{That is to say, for each neighborhood }$V$\textit{\ of }$0$ \textit{%
for the }$E$\textit{-topology, there exists }$\lambda_{0}\in\Lambda $\textit{%
\ such that }$\lambda\prec\lambda_{0}$ $\Longrightarrow$ $a_{\lambda
}(u_{\lambda}-v_{\lambda})\in V$. We write 
\begin{equation*}
u\overset{a}{\underset{E(\Omega)}{\sim}}v.
\end{equation*}
\end{definition}

\begin{remark}
\textit{We can also define an association process between }$u=\left[
u_{\lambda }\right] $\textit{\ and }$T\in \mathcal{E}(\Omega )$\textit{\ by
writing simply} 
\begin{equation*}
u\sim T\text{ }\Longleftrightarrow \text{ }\lim\limits_{E(\Omega ),\Lambda
}u_{\lambda }=T.
\end{equation*}%
Taking\ $E=\mathcal{D}^{\prime }$, $\mathcal{E}=\mathrm{C}^{\infty }$, $%
\Lambda =$ $(0,1],$ we recover the association process defined in the
literature ( \cite{Col1}, \cite{Col2}).
\end{remark}

\subsection{Algebraic framework for our problem}

Set $\mathcal{E}=\mathrm{C}^{\infty }$, $X=\mathbb{R}^{d}$ for $d=1,2,$ $E=%
\mathcal{D}^{\prime }$ and $\Lambda $ a set of indices, $\lambda \in \Lambda 
$. For any open set $\Omega $, in $\mathbb{R}^{d}$, $\mathcal{E}(\Omega )$
is endowed with the $\mathcal{P}(\Omega )$ topology of uniform convergence
of all derivatives on compact subsets of $\Omega $. This topology may be
defined by the family of the seminorms 
\begin{equation*}
P_{K,l}(u_{\lambda })=\sup_{\left\vert \alpha \right\vert \leq l}P_{K,\alpha
}(u_{\lambda })\text{ \ with }P_{K,\alpha }(u_{\lambda })=\sup_{x\in
K}\left\vert D^{\alpha }u_{\lambda }(x)\right\vert \text{, \ }K\Subset \Omega
\end{equation*}%
where the notation $\ K\Subset \mathbb{R}^{2}$ means that $K$ is a compact
subset of $\mathbb{R}^{2}$ and $D^{\alpha }=\dfrac{\partial ^{\alpha
_{1}+...+\alpha _{d}}}{\partial z_{1}^{\alpha _{1}}...\partial z_{d}^{\alpha
_{d}}}$ for $z=(z_{1},\ldots ,z_{d})\in \Omega $, $l\in \mathbb{N}$, $\alpha
=(\alpha _{1},...,\alpha _{d})\in \mathbb{N}^{d}$. Let $A$ be a subring of
the ring $\mathbb{R}^{\Lambda }$ of family of reals with the usual laws. We
consider a solid ideal$\ I_{A}$ of $A$. Then we have%
\begin{align*}
\mathcal{X}(\Omega )& =\{\left( u_{\lambda }\right) _{\lambda }\in \left[ 
\mathrm{C}^{\infty }(\Omega )\right] ^{\Lambda }:\forall K\Subset \Omega 
\text{, }\forall l\in \mathbb{N}\text{, }\left( P_{K,l}(u_{\lambda })\right)
_{\lambda }\in \left\vert A\right\vert \}\text{,} \\
\mathcal{N}(\Omega )& =\{\left( u_{\lambda }\right) _{\lambda }\in \left[ 
\mathrm{C}^{\infty }(\Omega )\right] ^{\Lambda }:\forall K\Subset \Omega 
\text{, }\forall l\in \mathbb{N}\text{, }\left( P_{K,l}(u_{\lambda })\right)
_{\lambda }\in \left\vert I_{A}\right\vert \}\text{,}\\
\mathcal{A}(\Omega )& =\mathcal{X}(\Omega )/\mathcal{N}(\Omega )\text{.}
\end{align*}%
The generalized derivation $D^{\alpha }:u(=\left[ u_{\varepsilon }\right]
)\mapsto D^{\alpha }u=\left[ D^{\alpha }u_{\varepsilon }\right] $ provides $%
\mathcal{A}(\Omega )$ with a differential algebraic structure (cf \cite{DelScar}).

\begin{example}
Set $\Lambda =(0,1]$. Consider 
\begin{equation*}
A=\mathbb{R}_{M}^{\Lambda }=\left\{ (m_{\lambda })_{\lambda }\in \mathbb{R}%
^{\Lambda }:\exists p\in \mathbb{R}_{+}^{\ast },\text{\ }\exists C\in 
\mathbb{R}_{+}^{\ast },\text{\ }\exists \mu \in \text{\ }(0,1],\text{\ }%
\forall \lambda \in \text{\ }(0,\mu ],\text{\ }\left\vert m_{\lambda
}\right\vert \leq C\lambda ^{-p}\right\}
\end{equation*}%
and the ideal 
\begin{equation*}
I_{A}=\left\{ (m_{\lambda })_{\lambda }\in \mathbb{R}^{\Lambda }:\forall
q\in \mathbb{R}_{+}^{\ast },\text{\ }\exists D\in \mathbb{R}_{+}^{\ast },%
\text{~}\exists \mu \in \text{\ }(0,1],\text{\ }\forall \lambda \in \text{\ }%
(0,\mu ],\text{\ }\left\vert m_{\varepsilon }\right\vert \leq D\lambda
^{q}\right\} \text{.}
\end{equation*}%
In this case we denote $\mathcal{X}^{s}(\Omega )=\mathcal{X}(\Omega )$ and $%
\mathcal{N}^{s}(\Omega )=\mathcal{N}(\Omega )$. The sheaf of factor algebras 
$\mathcal{G}^{s}\left( \mathbb{\cdot }\right) =\mathcal{X}^{s}(\cdot )/%
\mathcal{N}^{s}(\cdot )$ is called the sheaf of simplified Colombeau
algebras. $\mathcal{G}^{s}\left( \mathbb{R}^{d}\right) $ is the simplified
Colombeau algebra of generalized functions \cite{Col1}, \cite{Col2}.
\end{example}

We have the analogue of theorem 1.2.3. of \cite{GKOS} for $(\mathcal{C},%
\mathcal{E},\mathcal{P})$-algebras. We suppose here that $\Lambda $ is left
filtering and give this proposition for $\mathcal{A}\left( \mathbb{R}%
^{2}\right) $, although it is valid in more general situations.

\begin{proposition}
\label{0-estimate}Let $B$ be the set introduced in Definition \ref%
{Overgenerated} and assume that there exists $\left(
a_{\lambda}\right) _{\lambda}\in B$ with $\displaystyle \lim_{\lambda \to 0}a_{\lambda}=0$.\ C%
\textit{onsider }$\left( u_{\lambda}\right) _{\lambda}\in \mathcal{X}(%
\mathbb{R}^{2})$ such that%
\begin{equation*}
\forall K\Subset\mathbb{R}^{2},\ \ \left( P_{K,0}\left( u_{\lambda}\right)
\right) _{\lambda}\in\left\vert I_{A}\right\vert .
\end{equation*}
\textit{Then }$\left( u_{\lambda}\right) _{\lambda}\in\mathcal{N}(\mathbb{R}%
^{2})$.
\end{proposition}

We refer the reader to \cite{JAMADEV2}\ and \cite{Del2} for a detailed proof.

\begin{definition}
Tempered generalized functions, \cite{GKOS}, \cite{Scar}, \cite{Scar2}. For $%
f\in \mathrm{C}^{\infty }(\mathbb{R}^{n})$, $r\in \mathbb{Z}$ and $m\in 
\mathbb{N}$, we put%
\begin{equation*}
\mu _{r,m}(f)=\sup_{x\in \mathbb{R}^{n},\left\vert \alpha \right\vert \leq
m}(1+\left\vert x\right\vert )^{r}\left\vert \mathcal{D}^{\alpha
}f(x)\right\vert .
\end{equation*}%
The space of functions with slow growth \ is%
\begin{equation*}
\mathcal{O}_{M}(\mathbb{R}^{n})=\left\{ f\in \mathrm{C}^{\infty }(\mathbb{R}%
^{n}):\forall m\in \mathbb{N},\exists q\in \mathbb{N},\mu _{-q,m}(f)<+\infty
\right\} .
\end{equation*}
\end{definition}

\begin{definition}
We put%
\begin{align*}
\mathcal{X}_{\tau }\left( \mathbb{R}^{n}\right) & =\{\left( f_{\varepsilon
}\right) _{\varepsilon }\in \mathcal{O}_{M}(\mathbb{R}^{n})^{(0,1]}:\forall
m\in \mathbb{N}\text{, }\exists q\in \mathbb{N},\exists N\in \mathbb{N},\mu
_{-q,m}(f_{\varepsilon })=O(\varepsilon ^{-N})\text{ \ }(\varepsilon
\rightarrow 0)\}\text{,} \\
\mathcal{N}_{\tau }\left( \mathbb{R}^{n}\right) & =\{\left( f_{\varepsilon
}\right) _{\varepsilon }\in \mathcal{O}_{M}(\mathbb{R}^{n})^{(0,1]}:\forall
m\in \mathbb{N}\text{, }\exists q\in \mathbb{N},\forall p\in \mathbb{N},\mu
_{-q,m}(f_{\varepsilon })=O(\varepsilon ^{p})\text{ \ }(\varepsilon
\rightarrow 0)\}.
\end{align*}%
$\mathcal{X}_{\tau }\left( \mathbb{R}^{n}\right) $ is a subalgebra of $%
\mathcal{O}_{M}(\mathbb{R}^{n})^{(0,1]}$ and $\mathcal{N}_{\tau }\left( 
\mathbb{R}^{n}\right) $ an ideal of $\mathcal{X}_{\tau }\left( \mathbb{R}%
^{n}\right)$. The algebra $\mathcal{G}_{\tau }\left( \mathbb{R}^{n}\right) =%
\mathcal{X}_{\tau }\left( \mathbb{R}^{n}\right) /\mathcal{N}_{\tau }\left( 
\mathbb{R}^{n}\right) $ is called the algebra of tempered generalized
functions. The generalized derivation $\mathcal{D}^{\alpha }:u=\left[
u_{\varepsilon }\right] \mapsto \mathcal{D}^{\alpha }u=\left[ \mathcal{D}%
^{\alpha }u_{\varepsilon }\right] $ provides $\mathcal{G}_{\tau }\left( 
\mathbb{R}^{n}\right) $ with a differential a$\lg $ebraic structure.
\end{definition}

If $u$ is a generalized function of the variable $x\in \mathbb{R}^{2}$ and $%
F\in \mathrm{C}^{\infty }(\mathbb{R}^{3},\mathbb{R)}$, we extend the
notation $F(\cdot ,\cdot ,u)$ in the following way:

\begin{definition}
\textit{Let }$\Omega $ \textit{be an open subset of} $\mathbb{R}^{2}$ 
\textit{and} $F\in \mathrm{C}^{\infty }(\Omega \times \mathbb{R},\mathbb{R})$%
\textit{. We say that the algebra }$\mathcal{A}\left( \Omega \right) $\emph{%
\ is stable under}\textit{\ }$F$\textit{\ if the following two conditions
are satisfied:}

\begin{itemize}
\item \textit{For each }$K\Subset \mathbb{R}^{2}$\textit{, }$l\in \mathbb{N}$
and\textit{\ }$\left( u_{\varepsilon }\right) _{\varepsilon }\in \mathcal{X}%
(\Omega )$\textit{, there is a positive finite sequence }$C_{0}$\textit{%
,..., }$C_{l}$\textit{, such that}%
\begin{equation*}
P_{K,l}(F(\cdot ,\cdot ,u_{\varepsilon }))\leq
\dsum\limits_{i=0}^{l}C_{i}\left( P_{K,l}(u_{\varepsilon })\right) ^{i}.
\end{equation*}

\item \textit{For each }$K\Subset \mathbb{R}^{2}$\textit{, }$l\in \mathbb{N}$%
\textit{, }$\left( v_{\varepsilon }\right) _{\varepsilon }$ and $\left(
u_{\varepsilon }\right) _{\varepsilon }\in \mathcal{X}(\Omega )$\textit{,
there is a positive finite sequence }$D_{1}$\textit{,..., }$D_{l}$\textit{\
such that}
\end{itemize}
\end{definition}

\begin{equation*}
P_{K,l}(F(\cdot ,\cdot ,v_{\varepsilon })-F(\cdot ,\cdot ,u_{\varepsilon
})\leq \dsum\limits_{j=1}^{l}D_{j}\left( P_{K,l}(v_{\varepsilon
}-u_{\varepsilon })\right) ^{j}.
\end{equation*}

\begin{remark}
\textit{If }$\mathcal{A}\left( \Omega \right) $\textit{\ is stable under }$F$
\textit{then, for all }$\left( u_{\varepsilon }\right) _{\varepsilon }\in 
\mathcal{X}(\Omega )$\textit{\ and }$\left( i_{\varepsilon }\right)
_{\varepsilon }\in \mathcal{N}(\Omega )$\textit{, we have }$\left( F(\cdot
,\cdot ,u_{\varepsilon })\right) _{\varepsilon }\in \mathcal{X}(\Omega )$; $%
\left( F(\cdot ,\cdot ,u_{\varepsilon }+i_{\varepsilon })-F(\cdot ,\cdot
,u_{\varepsilon })\right) _{\varepsilon }\in \mathcal{N}(\Omega ).$
\end{remark}

\subsubsection{Generalized operator associated to a stability property}

If $\mathcal{A}\left( \mathbb{R}^{2}\right) $ if stable under $F$, for $u=%
\left[ u_{\varepsilon }\right] \in \mathcal{A}\left( \mathbb{R}^{2}\right) $%
, $\left[ F(.,.,u_{\varepsilon })\right] $ is a well defined element of $%
\mathcal{A}\left( \mathbb{R}^{2}\right) $ (i.e. not depending on the
representative $\left( u_{\varepsilon }\right) _{\varepsilon }$ of $u$).
This leads to the following:

\begin{definition}
\label{stability}If $\mathcal{A}\left( \mathbb{R}^{2}\right) $ if stable
under $F$, the operator 
\begin{equation*}
\mathcal{F}:\mathcal{A}\left( \mathbb{R}^{2}\right) \rightarrow \mathcal{A}%
\left( \mathbb{R}^{2}\right) ,\ \ u=\left[ u_{\varepsilon }\right] \mapsto %
\left[ F(.,.,u_{\varepsilon })\right]
\end{equation*}%
is called the \emph{generalized operator associated to }$F$. See \cite%
{JAMADEV2}.
\end{definition}

\subsubsection{Generalized restriction mappings}

Set $\left( f_{\varepsilon }\right) _{\varepsilon }$ be a family of
functions in $\mathrm{C}^{\infty }\left( \mathbb{R}\right) $. For each $g\in 
\mathrm{C}^{\infty }\left( \mathbb{R}^{2}\right) $ set 
\begin{equation*}
R_{\varepsilon }\left( g\right) :\mathrm{C}^{\infty }\left( \mathbb{R}%
\right) \rightarrow \mathrm{C}^{\infty }\left( \mathbb{R}\right) ,\
f_{\varepsilon }\mapsto \left( x\mapsto g(x,f_{\varepsilon }(x))\right) .
\end{equation*}%
The family $\left( R_{\varepsilon }\right) _{\varepsilon }$ map $\left( 
\mathrm{C}^{\infty }\left( \mathbb{R}^{2}\right) \right) ^{\Lambda }$ into $%
\left( \mathrm{C}^{\infty }\left( \mathbb{R}\right) \right) ^{\Lambda }$.

\begin{definition}
The family of smooth function $\left( f_{\varepsilon }\right) _{\varepsilon
} $ is \emph{compatible with second side restriction} if%
\begin{eqnarray*}
\forall \left( u_{\varepsilon }\right) _{\varepsilon } &\in &\mathcal{X}((%
\mathbb{R}^{2}),\ \ \left( u_{\varepsilon }\left( \cdot ,f_{\varepsilon
}(\cdot )\right) \right) _{\varepsilon }\in \mathcal{X}(\mathbb{R})\ ; \\
\forall \left( i_{\varepsilon }\right) _{\varepsilon } &\in &\mathcal{N}(%
\mathbb{R}^{2}),\ \ \left( i_{\varepsilon }\left( \cdot ,f_{\varepsilon
}(\cdot )\right) \right) _{\varepsilon }\in \mathcal{N}(\mathbb{R}).
\end{eqnarray*}
\end{definition}

Clearly, if $u=\left[ u_{\varepsilon }\right] \in \mathcal{A}(\mathbb{R}%
^{2}) $ then $\left[ u_{\varepsilon }\left( \cdot ,f_{\varepsilon }(\cdot
)\right) \right] $ is a well defined element of $\mathcal{A}(\mathbb{R)}$
(i.e. not depending on the representative of $u$.) This leads to the
following:

\begin{definition}
\label{restriction}If the family of smooth function $\left( f_{\varepsilon
}\right) _{\varepsilon }$ is compatible with second side restriction, the
mapping%
\begin{equation*}
\mathcal{R}:\mathcal{A}\left( \mathbb{R}^{2}\right) \rightarrow \mathcal{A}%
\left( \mathbb{R}\right) ,\ \ \ u=\left[ u_{\varepsilon }\right] \mapsto %
\left[ u_{\varepsilon }\left( \cdot ,f_{\varepsilon }(\cdot )\right) \right]
=\left[ R_{\varepsilon }\left( u_{\varepsilon }\right) \right]
\end{equation*}%
is called the\emph{\ generalized second side restriction mapping} associated
to the family $\left( f_{\varepsilon }\right) _{\varepsilon }$.
\end{definition}

\begin{remark}
The previous process generalizes the standard one defining the restriction
of the generalized function $u=\left[ u_{\varepsilon }\right] \in \mathcal{A}%
\left( \mathbb{R}^{2}\right) $ to the manifold $\left\{ y=f\left( x\right)
\right\} $ obtained when taking $f_{\varepsilon }=f$ for each $\varepsilon
\in \Lambda $.
\end{remark}

First let us state a useful definition used throughout this article:

\begin{definition}
\cite{GKOS} Let $(f_{\varepsilon })_{\varepsilon }$ be a family of $%
C^{\infty }(\mathbb{R}^{n})$ functions. This family is c-bounded if for all
compact set $K\subset \mathbb{R}^{n}$ it exists another compact set $%
L\subset \mathbb{R}^{n}$ such that $f_{\varepsilon }(K)\subset L$ for all $%
\varepsilon $ ($L$ is independent of $\varepsilon $).
\end{definition}

\begin{proposition}
Assume that:\newline
$\left( i\right) ~$ For each $K\Subset \mathbb{R}$, it exists $K^{\prime
}\Subset \mathbb{R}$ such that, for all $\varepsilon \in \Lambda $, $%
f_{\varepsilon }(K)\subset K^{\prime },$\newline
$\left( ii\right) ~\left( f_{\varepsilon }\right) _{\varepsilon }$ belongs
to $\mathcal{X}(\mathbb{R})$.\newline
Then the family $\left( f_{\varepsilon }\right) _{\varepsilon }$ is
compatible with restriction.
\end{proposition}

\begin{proof}
Take $\left( u_{\varepsilon }\right) _{\varepsilon }$ (resp. $\left(
i_{\varepsilon }\right) _{\varepsilon }$) in $\mathcal{X}(\mathbb{R}^{2})$
(resp. $\mathcal{N}(\mathbb{R}^{2})$) and set $v_{\varepsilon }\left(
x\right) =u_{\varepsilon }( x,f_{\varepsilon }(x)) $. From\ $%
\left( i\right) $ we have 
\begin{align*}
p_{K,0}\left( v_{\varepsilon }\right) & \leq p_{K\times K^{\prime },0}\left(
u_{\varepsilon }\right)  ,\\
p_{K,1}\left( v_{\varepsilon }\right) & \leq p_{K\times K^{\prime },\left(
1,0\right) }\left( u_{\varepsilon }\right) +p_{K\times K^{\prime },\left(
0,1\right) }\left( u_{\varepsilon }\right) p_{K,1}\left( f_{\varepsilon
}\right) .
\end{align*}%
By induction we can see that for each\textit{\ }$K\Subset \mathbb{R}$\textit{%
, }and each\textit{\ }$l\in \mathbb{N}$, $p_{K,l}\left( v_{\varepsilon
}\right) $ is estimated by sums or products of terms like $p_{K\times
K^{\prime },\left( n,m\right) }\left( u_{\varepsilon }\right) $ for $n+m\leq
l$, or $p_{K,k}\left( f_{\varepsilon }\right) $ for $k\leq l$. Then, from $%
\left( ii\right) $, $p_{K,l}\left( v_{\varepsilon }\right) $ is in $%
\left\vert A\right\vert $. Similarly, setting $j_{\varepsilon }\left(
x\right) =i_{\varepsilon }\left( x,f_{\varepsilon }(x\right) $ leads to $%
p_{K,l}\left( j_{\varepsilon }\right) \in \left\vert I_{A}\right\vert $.
Then $\left( u_{\varepsilon }\left( .,f_{\varepsilon }(.\right) \right)
_{\varepsilon }$ (resp. $i_{\varepsilon }\left( .,f_{\varepsilon }(.\right) $%
) belongs to $\mathcal{X}(\mathbb{R})$ (resp. $\mathcal{N}(\mathbb{R})$).
\end{proof}

\section{Existence of solutions for a characteristic Cauchy problem in $(%
\mathcal{C},\mathcal{E},\mathcal{P})$-algebras}

\bigskip \label{sect_existence}

We will use the notations found in \cite{Dev1} and \cite{Dev2}.

\subsection{A generalized differential problem associated to the ill posed
classical one}

Our goal is to give a meaning to the characteristic Cauchy problem formally
written as%
\begin{equation*}
\left( P_{form}\right) \left\{ 
\begin{array}{l}
\dfrac{\partial ^{2}u}{\partial x\partial y}=F(\cdot ,\cdot ,u), \\ 
\left. u\right\vert _{\gamma }=\varphi , \\ 
\left. \dfrac{\partial u}{\partial y}\right\vert _{\gamma }=\psi ,%
\end{array}%
\right.
\end{equation*}%
where the data $\psi $, $\varphi $ are smooth functions given on a
characteristic \textrm{C}$^{\infty }$ manifold $\gamma $ supposed to be a
curve whose equation is $y=f(x)$. $F$ is smooth in its arguments.

We don't have a classical surrounding in which we can pose (and a fortiori
solve) the problem. In the sequel, by means of regularizing processes we
will define an associated problem to $\left( P_{form}\right) $.%
\begin{equation*}
\left( P_{gen}\right) \left\{ 
\begin{array}{c}
\dfrac{\partial ^{2}u}{\partial x\partial y}=\mathcal{F}(u), \\ 
\mathcal{R}\left( u\right) =\varphi , \\ 
\mathcal{R}\left( \dfrac{\partial u}{\partial y}\right) =\psi%
\end{array}%
\right.
\end{equation*}%
where $u$ is searched in some convenient algebra $\mathcal{A}\left( \mathbb{R%
}^{2}\right) $, $\mathcal{F}$, $\mathcal{R}$ are defined as
previously.

The idea is then to approach this Cauchy problem by a family of
non-characteristic ones by replacing the characteristic curve $\gamma $ by a
family of smooth non-characteristic curves $\gamma _{\varepsilon }$ whose
equation is $y=f_{\varepsilon }(x)$. Moreover $\gamma _{\varepsilon }$ is
diffeomorphic to $y=0$, which is a consequence of the following assumption%
\begin{equation*}
(H):\left\{ 
\begin{array}{l}
F\in \mathrm{C}^{\infty }(\mathbb{R}^{3},\mathbb{R)}, \\ 
\forall K\Subset \mathbb{R}^{2}\text{, }\sup_{\left( x,y\right) \in K;z\in 
\mathbb{R}}\text{ }\left\vert \partial _{z}F(x,y,z)\right\vert <\infty , \\ 
f\text{ }_{\varepsilon }\text{ is defined and strictly increasing on }%
\mathbb{R}\text{ with image }\mathbb{R}, \\ 
\forall x\in \mathbb{R},f_{\varepsilon }^{\prime }(x)\neq 0\text{.}%
\end{array}%
\right.
\end{equation*}

In terms of representatives, and thanks to the stability and restriction
hypothesis, solving $\left( P_{gen}\right) $ amounts to find a family $%
\left( u_{\varepsilon }\right) _{\varepsilon }\in \mathcal{X}(\mathbb{R}^{2})
$ such that 
\begin{equation*}
\left\{ 
\begin{array}{c}
\dfrac{\partial ^{2}u_{\varepsilon }}{\partial x\partial y}%
(x,y)-F_{\varepsilon }(x,y,u_{\varepsilon }\left( x,y\right)
)=i_{\varepsilon }\left( x,y\right) , \\ 
u_{\varepsilon }\left( x,f_{\varepsilon }(x)\right) -\varphi \left( x\right)
=j_{\varepsilon }\left( x\right) , \\ 
\dfrac{\partial u_{\varepsilon }}{\partial y}\left( x,f_{\varepsilon
}(x)\right) -\psi \left( x\right) =l_{\varepsilon }\left( x\right) ,%
\end{array}%
\right. 
\end{equation*}%
where $\left( i_{\varepsilon }\right) _{\varepsilon }\in \mathcal{N}\left( 
\mathbb{R}^{2}\right) $, $\left( j_{\varepsilon }\right) _{\varepsilon }$, $%
\left( l_{\varepsilon }\right) _{\varepsilon }\in \mathcal{N}\left( \mathbb{R%
}\right) $. Suppose we can find $u_{\varepsilon }\in \mathrm{C}^{\infty
}\left( \mathbb{R}^{2}\right) $ verifying

\begin{equation*}
\left\{ 
\begin{array}{c}
\dfrac{\partial ^{2}u_{\varepsilon }}{\partial x\partial y}%
(x,y)=F_{\varepsilon }(x,y,u_{\varepsilon }\left( x,y\right) ), \\ 
u_{\varepsilon }\left( x,f_{\varepsilon }(x)\right) =\varphi \left( x\right)
, \\ 
\dfrac{\partial u_{\varepsilon }}{\partial y}\left( x,f_{\varepsilon
}(x)\right) =\psi \left( x\right) ,%
\end{array}%
\right. 
\end{equation*}%
then, if we can prove that $\left( u_{\varepsilon }\right) _{\varepsilon
}\in \mathcal{X}(\mathbb{R}^{2})$, $u=\left[ u_{\varepsilon }\right] $ is a
solution of $\left( P_{gen}\right) $. If $v=\left[ v_{\varepsilon }\right] $
is another solution of $\left( P_{gen}\right) $ obtain by replacing $\gamma $
by another family of smooth non-characteristic curves $\gamma _{\varepsilon
}^{\prime }$ whose equation is $y=g_{\varepsilon }(x)$, this implies%
\begin{equation*}
\left\{ 
\begin{array}{c}
\dfrac{\partial ^{2}(v_{\varepsilon }-u_{\varepsilon })}{\partial x\partial y%
}(x,y)-(F_{\varepsilon }(x,y,v_{\varepsilon }\left( x,y\right)
)-F_{\varepsilon }(x,y,u_{\varepsilon }\left( x,y\right) ))=a_{\varepsilon
}\left( x,y\right) , \\ 
v_{\varepsilon }\left( x,g_{\varepsilon }(x)\right) -u_{\varepsilon }\left(
x,f_{\varepsilon }(x)\right) =b_{\varepsilon }\left( x\right) , \\ 
\dfrac{\partial v_{\varepsilon }}{\partial y}\left( x,g_{\varepsilon
}(x)\right) -\dfrac{\partial u_{\varepsilon }}{\partial y}\left(
x,f_{\varepsilon }(x)\right) =c_{\varepsilon }\left( x\right) ,%
\end{array}%
\right. 
\end{equation*}%
where $\left( a_{\varepsilon }\right) _{\varepsilon }\in \mathcal{N}\left( 
\mathbb{R}^{2}\right) $ and $\left( b_{\varepsilon }\right) _{\varepsilon }$%
, $\left( c_{\varepsilon }\right) _{\varepsilon }\in \mathcal{N}\left( 
\mathbb{R}\right) $. We have to prove that $\left( v_{\varepsilon
}-u_{\varepsilon }\right) _{\varepsilon }\in \mathcal{N}(\mathbb{R}^{2})$ if
we intend to prove that the solution of $\left( P_{gen}\right) $
in the algebra $\mathcal{A}\left( \mathbb{R}^{2}\right) $ does not depend on the representative of the class $[f_\e]$.

\subsection{Notations, assumptions and results}

Using the results of \cite{Dev2}, with assumption $\left( H\right) $ we have
a unique smooth solution of $\left( P_{\varepsilon }\right) $ which
satisfies the following integral equation%
\begin{equation*}
u_{\varepsilon }(x,y)=u_{0,\varepsilon
}(x,y)-\diint\nolimits_{D(x,y,f_{\varepsilon })}F(\xi ,\eta ,u_{\varepsilon
}(\xi ,\eta ))d\xi d\eta ,
\end{equation*}%
with $u_{0,\varepsilon }(x,y)=\varphi (x)-\chi _{\varepsilon }\left(
f_{\varepsilon }\left( x\right) \right) +\chi _{\varepsilon }\left( y\right) 
$, where $\chi _{\varepsilon }$ is a primitive of $\psi \circ f_{\varepsilon
}^{-1}$ and 
\begin{equation*}
D(x,y,f_{\varepsilon })=\left\{ 
\begin{array}{c}
\left\{ (\xi ,\eta )/x\leq \xi \leq f_{\varepsilon }^{-1}\left( y\right)
,f_{\varepsilon }(\xi )\leq \eta \leq y\right\} \text{, if }y\geq
f_{\varepsilon }\left( x\right) \text{,} \\ 
\left\{ (\xi ,\eta )/f_{\varepsilon }^{-1}\left( y\right) \leq \xi \leq
x,y\leq \eta \leq f_{\varepsilon }(\xi )\right\} \text{, if }y\leq
f_{\varepsilon }\left( x\right) .%
\end{array}%
\right.
\end{equation*}%
\textbf{Remarks, notations and hypothesis.}

\noindent Each compact $K\Subset \mathbb{R}^{2}$ is contained in some
product $\left[ -a,a\right] \times \left[ -b,b\right] $. We define 
\begin{equation}
\begin{cases}
\beta _{K,\varepsilon }=\max (a,f_{\varepsilon }^{-1}(b))\text{ and }\alpha
_{K,\varepsilon }=\min (-a,f_{\varepsilon }^{-1}(-b)), \\ 
a_{K,\varepsilon }=2\max (\beta _{K,\varepsilon },\left\vert \alpha
_{K,\varepsilon }\right\vert ), \\ 
K_{\varepsilon }=K_{1\varepsilon }\times K_{2}\text{ with }K_{1\varepsilon }=%
\left[ -a_{K,\varepsilon }/2,a_{K,\varepsilon }/2\right] \text{ and }K_{2}=%
\left[ -b,b\right] =\left[ -c/2,c/2\right] .%
\end{cases}
\label{4}
\end{equation}

By construction we have $K\subset K_{\varepsilon }$.

We also make the following assumptions to generate a convenient $(\mathcal{C}%
,\mathcal{E},\mathcal{P})$-algebra adapted to our problem:%
\begin{equation}
\left( H_{1}\right) \left\{ 
\begin{array}{c}
\forall \left( \varepsilon ,\eta \right) \in (0,1]^{2},\forall K\Subset 
\mathbb{R}^{2},\forall l\in \mathbb{N},\exists \mu _{K,l}>0,\exists
M_{\varepsilon }>0, \\ 
\underset{(t,x,z)\in K_{\varepsilon }\times \mathbb{R},\text{ }\left\vert
\alpha \right\vert \leq l}{\sup }\left\vert D^{\alpha }F(t,x,z)\right\vert
=M_{K,\varepsilon ,l}\leq \mu _{K,l}M_{\varepsilon }. \\ 
\forall \varepsilon \in (0,1],\forall K\Subset \mathbb{R}^{2},\exists \nu
_{K}>0,\exists a_{\varepsilon }>0,a_{K,\varepsilon }\leq \nu
_{K}a_{\varepsilon }.%
\end{array}%
\right.
\end{equation}%
Particularly, we set 
\begin{equation*}
m_{K,\varepsilon }=\underset{(t,x)\in K_{\varepsilon };\text{ }z\in \mathbb{R%
}}{\sup }\left\vert \dfrac{\partial }{\partial z}F(t,x,z)\right\vert \leq
\mu _{K,1}M_{\varepsilon }.
\end{equation*}

\begin{align*}
\left( H_{2}\right) & 
\begin{cases}
\exists \left( r_{\varepsilon }\right) _{\varepsilon }\in \mathbb{R}_{\ast
}^{(0,1]}\text{ such that }\forall K_{2}\Subset \mathbb{R},\forall \alpha
_{2}\in \mathbb{N},\exists D_{2}=D_{K_{2},\alpha _{2},\varepsilon }\in 
\mathbb{R}_{+}^{\ast },\exists p\in \mathbb{N}, \\ 
\max \left[ \underset{K_{2}}{\sup }\left\vert D^{\alpha _{2}}\varphi
(f_{\varepsilon }^{-1}\left( y\right) )\right\vert ,\underset{K_{2}}{\sup }%
\left\vert D^{\alpha _{2}}\chi _{\varepsilon }\left( y\right) \right\vert %
\right] \leq \frac{D_{2}}{\left( r_{\varepsilon }\right) ^{p}}.%
\end{cases}
\\
\left( H_{3}\right) & 
\begin{cases}
\mathcal{C}=A/I_{A}\text{ is overgenerated by the following elements of }%
\mathbb{R}_{\ast }^{(0,1]} \\ 
\left( \varepsilon \right) _{\varepsilon },\left( r_{\varepsilon }\right)
_{\varepsilon },\left( M_{\varepsilon }\right) _{\varepsilon },(\exp
M_{\varepsilon }a_{\varepsilon }).%
\end{cases}
\\
\left( H_{4}\right) & 
\begin{cases}
\mathcal{A}\left( \mathbb{R}^{2}\right) =\mathcal{X}(\mathbb{R}^{2})/%
\mathcal{N}(\mathbb{R}^{2})\text{ is built on }\mathcal{C}\text{ with} \\ 
\text{ }(\mathcal{E},\mathcal{P})=\left( \mathrm{C}^{\infty }(\mathbb{R}%
^{2}),\left( P_{K,l}\right) _{K\Subset \mathbb{R}^{2},l\in \mathbb{N}%
}\right)  \\ 
\text{and }\mathcal{A}\left( \mathbb{R}^{2}\right) \text{ is stable under }F%
\text{ relatively to }\mathcal{C}.%
\end{cases}
\\
\left( H_{5}\right) & \text{ \ \ \ \ \ \ }f_{\varepsilon }\in \mathrm{C}%
^{\infty }(\mathbb{R})\text{, }f_{\varepsilon }\text{ strictly increasing, }%
f_{\varepsilon }(\mathbb{R})=\mathbb{R}\text{ and }\phi,\psi \in O_M(\R) \\
\left( H_{6}\right) & \text{ \ \ \ \ \ }(f_{\varepsilon })_{\varepsilon
},(f_{\varepsilon }^{-1})_{\varepsilon }\in \mathcal{X}_{\tau }\left( 
\mathbb{R}\right) \text{, }(f_{\varepsilon })_{\varepsilon }\text{ is
c-bounded}\text{ and }\underset{\varepsilon \underset{\mathcal{D}^{\prime }(%
\mathbb{R})}{\rightarrow }0}{\lim }f_{\varepsilon }=f.
\end{align*}

\begin{lemma}
\label{relations} We have the following relations: 
\begin{align}
& \left( \left\vert \alpha _{K,\varepsilon }\right\vert \right)
_{\varepsilon },\left( \left\vert \beta _{K,\varepsilon }\right\vert \right)
_{\varepsilon },\left( \left\vert a_{K,\varepsilon }\right\vert \right)
_{\varepsilon }\in \left\vert A\right\vert ,  \label{A6} \\
& \forall \varepsilon ,\forall (x,y)\in K_{\varepsilon
},D(x,y,f_{\varepsilon })\subset K_{\varepsilon }.  \label{A7}
\end{align}
\end{lemma}

%
%
%
%
%
%
%
%
%
%
%
%
First $(f_{\varepsilon }^{-1})_{\varepsilon }\in $ $\mathcal{X_{\tau }}(%
\mathbb{R})$ so $\left( \left\vert \alpha _{K,\varepsilon }\right\vert
\right) _{\varepsilon },\left( \left\vert \beta _{K,\varepsilon }\right\vert
\right) _{\varepsilon }\in \left\vert A\right\vert $ and then obviously $%
\left( \left\vert a_{K,\varepsilon }\right\vert \right) _{\varepsilon }\in
\left\vert A\right\vert $. Next as $(f_{\varepsilon })_{\varepsilon }\in 
\mathcal{X}_{\tau }(\mathbb{R})$ we can find $p\in \mathbb{N}$ such that $%
\forall x,\varepsilon ,\left\vert f_{\varepsilon }(x)\right\vert \leq
\varepsilon ^{-p}(1+\left\vert x\right\vert )^{p}$ so we have 
\begin{equation*}
\left\vert \mu _{\varepsilon }\right\vert =\left\vert f_{\varepsilon
}(\alpha _{K,\varepsilon })\right\vert \leq \varepsilon ^{-p}(1+\left\vert
\alpha _{K,\varepsilon }\right\vert )^{p}\in \left\vert A\right\vert .
\end{equation*}

\begin{theorem}
\label{thm_solution_moderee} With the notations and the hypothesis of the
above paragraph, if $u_{\varepsilon }$ is the solution to the problem $%
\left( P_{\varepsilon }\right) $ then problem $\left( P_{gen}\right) $
admits $u=\left[ u_{\varepsilon }\right] _{\mathcal{A}(\mathbb{R}^{2})}$ as
solution.
\end{theorem}

\textbf{Proof.} \ \ \ We have: $u_{\varepsilon }(x,y)=u_{0,\varepsilon
}(x,y)-u_{1,\varepsilon }(x,y)$, where 
\begin{equation*}
u_{0,\varepsilon }(x,y)=\varphi (x)-\chi _{\varepsilon }\left(
f_{\varepsilon }\left( x\right) \right) +\chi _{\varepsilon }\left( y\right) 
\text{,}
\end{equation*}%
and 
\begin{equation*}
u_{1,\varepsilon }(x,y)={\displaystyle\iint\nolimits_{D(x,y,f_{\varepsilon
})}}F(\xi ,\eta ,u_{\varepsilon }(\xi ,\eta ))d\xi d\eta .
\end{equation*}%
We will actually prove that $(P_{K_{\varepsilon },n}(u_{\varepsilon
}))_{\varepsilon }\in \left\vert A\right\vert $.\newline
First we have $\chi _{\varepsilon }^{\prime }=\psi \circ f_{\varepsilon
}^{-1}$; as for $f_{\varepsilon }^{-1}(K_{2}])=K_{1\varepsilon }$ and as $%
\psi \in \mathcal{O}_{M}(\mathbb{R})$, $\left( \left\vert \alpha
_{K,\varepsilon }\right\vert \right) _{\varepsilon },\left( \left\vert \beta
_{K,\varepsilon }\right\vert \right) _{\varepsilon }\in \left\vert
A\right\vert $ then 
\begin{equation*}
\forall l\in \mathbb{N},\left( P_{K_{2},l}(\chi _{\varepsilon })\right)
_{\varepsilon }\in \left\vert A\right\vert .
\end{equation*}%
Moreover as $\varphi \in \mathcal{O}_{M}(\mathbb{R})$ we also have that 
\begin{equation*}
\forall l\in \mathbb{N},\left( P_{K_{1\varepsilon },l}(\varphi )\right)
_{\varepsilon }\in \left\vert A\right\vert
\end{equation*}%
and as $(\chi _{\varepsilon }\circ f_{\varepsilon })^{\prime
}=f_{\varepsilon }^{\prime }\psi $ and $(f_{\varepsilon }^{\prime
})_{\varepsilon }\in \mathcal{X}_{\tau }(\mathbb{R})$ we can conclude that 
\begin{equation*}
\forall l\in \mathbb{N},\left( P_{K_{\varepsilon },l}\left( u_{0,\varepsilon
}\right) \right) _{\varepsilon }\in \left\vert A\right\vert .
\end{equation*}%
We have the following equality 
\begin{equation*}
F(x,y,u_{\varepsilon }(x,y))-F(x,y,0)=u_{\varepsilon }(x,y)\int_{0}^{1}\frac{%
\partial {F}}{\partial {z}}(x,y,\theta u_{\varepsilon }(x,y))d\theta 
\end{equation*}%
so that 
\begin{equation*}
\left\vert u_{\varepsilon }(x,y)\right\vert \leq \left\vert \left\vert
u_{0,\varepsilon }\right\vert \right\vert _{K_{\varepsilon
}}+\iint_{D(x,y,f_{\varepsilon })}\left\vert u_{\varepsilon }(\xi ,\eta
)\int_{0}^{1}\frac{\partial {F}}{\partial {z}}(x,y,\theta u_{\varepsilon
}(x,y))d\theta +F(x,y,0)\right\vert d\xi d\eta .
\end{equation*}%
Now if we define ($A(K_{\varepsilon })$ denotes the area of $K_{\varepsilon }
$) 
\begin{equation*}
C_{\varepsilon }=A(K_{\varepsilon })\sup_{(x,y)\in D_{\varepsilon
}}\left\vert F(x,y,0)\right\vert +\left\vert \left\vert u_{0,\varepsilon
}\right\vert \right\vert _{K_{\varepsilon }}\in \left\vert A\right\vert ,
\end{equation*}%
we get (remembering relation (\ref{A7})) 
\begin{equation*}
\forall (x,y)\in K_{\varepsilon },\left\vert u_{\varepsilon
}(x,y)\right\vert \leq C_{\varepsilon }+\iint_{D(x,y,f_{\varepsilon })}\mu
_{K,1}M_{\varepsilon }\left\vert u_{\varepsilon }(\xi ,\eta )\right\vert
d\xi d\eta .
\end{equation*}%
We define $e_{\varepsilon }(y)=\sup_{x\in K_{1\varepsilon }}\left\vert
u_{\varepsilon }(x,y)\right\vert $ and then for all $y\in K_{2}$ we have 
\begin{equation*}
e_{\varepsilon }(y)\leq C_{\varepsilon }+\int_{-b}^{b}\mu
_{K,1}M_{\varepsilon }e_{\varepsilon }(\eta )d\eta ,
\end{equation*}%
so%
\begin{equation*}
e_{\varepsilon }(y)\leq C_{\varepsilon }\exp \left( 2b\mu
_{K,1}M_{\varepsilon }\right) \text{ using Gronwall's lemma}
\end{equation*}%
and then as $C_{\varepsilon }\exp \left( 2b\mu _{K,1}M_{\varepsilon }\right)
\in \left\vert A\right\vert $ we have proved that 
\begin{equation*}
(P_{K_{\varepsilon },0}(u_{\varepsilon }))_{\varepsilon }=(\left\vert
\left\vert u_{\varepsilon }\right\vert \right\vert _{K_{\varepsilon
}})_{\varepsilon }\in \left\vert A\right\vert .
\end{equation*}%
Let us assume now that we have $\left( P_{K_{\varepsilon },n}(u_{\varepsilon
})\right) _{\varepsilon }\in \left\vert A\right\vert $, for all $n\in 
\mathbb{N}$.\newline
By successive derivations, for $n\geq 1$ we obtain 
\begin{eqnarray}
\frac{\partial ^{n+1}u_{1,\varepsilon }}{\partial x^{n+1}}(x,y) &=&-{%
\displaystyle\sum\nolimits_{j=0}^{n-1}}C_{n}^{j}f_{\varepsilon }^{\left(
n-j\right) }(x)\frac{\partial ^{j}}{\partial x^{j}}F(x,f_{\varepsilon
}(x),\varphi (x))  \label{A8} \\
&&+\int_{f_{\varepsilon }(x)}^{y}\frac{\partial ^{n}}{\partial x^{n}}%
F(x,\eta ,u_{\varepsilon }(x,\eta ))d\eta .  \label{A9}
\end{eqnarray}%
First as $(\alpha _{\varepsilon })_{\varepsilon },(\beta _{\varepsilon
})_{\varepsilon }\in \left\vert A\right\vert $ and $\varphi \in \mathcal{O}%
_{M}(\mathbb{R})$ we have that $(P_{K_{1\varepsilon },j}(\varphi
))_{\varepsilon }\in \left\vert A\right\vert $ because for any $k\in \mathbb{%
N}$, we can find $p\in \mathbb{N}$ such that $\left\vert \varphi
^{(k)}(\alpha _{\varepsilon })\right\vert \leq (1+\left\vert \alpha
_{\varepsilon }\right\vert )^{p}$.\newline
Moreover $f_{\varepsilon }\in \mathcal{X}_{\tau }$ then for all $k$, we can
find $p\in \mathbb{N}$ such that 
\begin{equation*}
\forall \varepsilon ,\sup_{\mathbb{R}}(1+\left\vert x\right\vert
)^{-p}\left\vert f_{\varepsilon }^{(k)}(x)\right\vert \leq \varepsilon ^{-p},
\end{equation*}%
but then we have 
\begin{equation*}
\left\vert \left\vert f_{\varepsilon }^{(k)}\right\vert \right\vert
_{K_{1\varepsilon }}\leq \max \left\{ (1+\left\vert \alpha _{\varepsilon
}\right\vert )^{p},(1+\left\vert \beta _{\varepsilon }\right\vert
)^{p}\right\} \varepsilon ^{-p}\in \left\vert A\right\vert
\end{equation*}%
and as $(P_{K_{\varepsilon }\times \mathbb{R},j}(F))_{\varepsilon }\in
\left\vert A\right\vert $ this takes care of the first term in (\ref{A8}).
For the second term we compute first for $n=1$ 
\begin{eqnarray*}
&&\int_{f_{\varepsilon }(x)}^{y}\frac{\partial }{\partial x}F(x,\eta
,u_{\varepsilon }(x,\eta ))d\eta  \\
&=&\int_{f_{\varepsilon }(x)}^{y}\left( \frac{\partial F}{\partial x}(x,\eta
,u_{\varepsilon }(x,\eta ))+\frac{\partial F}{\partial z}(x,\eta
,u_{\varepsilon }(x,\eta ))\frac{\partial u_{\varepsilon }}{\partial x}%
(x,\eta )\right) d\eta ,
\end{eqnarray*}%
then we have 
\begin{equation*}
\left\vert \left\vert \int_{f_{\varepsilon }(x)}^{y}\frac{\partial }{%
\partial x}F(x,\eta ,u_{\varepsilon }(x,\eta ))d\eta \right\vert \right\vert
_{D_{\varepsilon }}\leq 2bP_{K_{\varepsilon }\times \mathbb{R},1}(F)\left(
1+P_{K_{\varepsilon },1}(u_{\varepsilon })\right) ,
\end{equation*}%
and then this is in $\left\vert A\right\vert $ because of the hypothesis for 
$F$. Now for $n=2$ we have similarly 
\begin{equation*}
\left\vert \left\vert \int_{f_{\varepsilon }(x)}^{y}\frac{\partial ^{2}}{%
\partial x^{2}}F(x,\eta ,u_{\varepsilon }(x,\eta ))d\eta \right\vert
\right\vert _{D_{\varepsilon }}\leq 2bP_{K_{\varepsilon }\times \mathbb{R}%
,2}(F)\left( P_{K_{\varepsilon },1}(u_{\varepsilon })^{2}+P_{K_{\varepsilon
},2}(u_{\varepsilon })+1\right) 
\end{equation*}%
and this is also in $\left\vert A\right\vert $ as the induction hypothesis
insures that $(P_{K_{\varepsilon },2}(u_{\varepsilon }))_{\varepsilon }\in
\left\vert A\right\vert $.\newline
For $n>2$ the same calculation leads to more terms involving higher
derivatives but each of them can be dealt with using the same kind of
arguments.\newline
Similarly we have 
\begin{eqnarray*}
\frac{\partial ^{n+1}u_{1,\varepsilon }}{\partial y^{n+1}}(x,y)
&=&-\sum\nolimits_{j=0}^{n-1}C_{n}^{j}\left( f_{\varepsilon }^{-1}\right)
^{\left( n-j\right) }(y)\frac{\partial ^{j}}{\partial y^{j}}F(f_{\varepsilon
}^{-1}(y),y,\varphi (f_{\varepsilon }^{-1}(y))) \\
&&-\int_{x}^{f_{\varepsilon }^{-1}(y)}\frac{\partial ^{n}}{\partial y^{n}}%
F(\xi ,y,u_{\varepsilon }(\xi ,y))d\xi .
\end{eqnarray*}%
As $f_{\varepsilon }^{-1}\in \mathcal{X}_{\tau }(\mathbb{R})$ and using the
same argument as previously the first terms have their $P_{K}$ norm in $%
\left\vert A\right\vert $.\newline
Moreover we have 
\begin{equation*}
\left\vert \left\vert \int_{x}^{f_{\varepsilon }^{-1}(y)}\frac{\partial ^{n}%
}{\partial y^{n}}F(\xi ,y,u_{\varepsilon }(\xi ,y))d\xi \right\vert
\right\vert _{K_{\varepsilon }}\leq a_{K,\varepsilon }\left\vert \left\vert 
\frac{\partial ^{n}}{\partial y^{n}}F(x,y,u_{\varepsilon }(x,y))\right\vert
\right\vert _{K_{\varepsilon }}.
\end{equation*}%
The same arguments apply here so that this term also has its $%
P_{K_{\varepsilon }}$ norm in $\left\vert A\right\vert $.\newline
The proof for the other partial derivatives can be done along the same
lines, thus finally we can conclude that $\left( P_{K_{\varepsilon
},n+1}(u_{\varepsilon })\right) _{\varepsilon }\in \left\vert A\right\vert $
which concludes the induction. 

It is to be observed that we prove more than we need for the existence
alone, but we will definitely need this stronger statement when proving that this solution only depends on the class $[f_\e]$ as we will have to compare different generalized solutions.

\subsection{Generalized solutions only depend on the class $[f_\e]$}

\label{sect_unicite}

\begin{theorem}
\label{thm_solution_unique} Under the same hypotheses as theorem \ref%
{thm_solution_moderee}, the generalized function $u$ represented by the
family $\left( u_{\varepsilon }\right) _{\varepsilon }$ of solutions to
Problems $\left( P_{\varepsilon }\right) $, does not depend on the choice of
the representative $\left( f_{\varepsilon }\right) _{\varepsilon }$ of the
class $f=\left[ f_{\varepsilon }\right] \in \mathcal{G}_{\tau }\left( 
\mathbb{R}\right) $.
\end{theorem}

\begin{remark}
\label{rem_cont_reciprocal} We need to consider tempered generalized
functions for $(f_\e)_\e$ as this counter-example shows:\newline
Take $f_{\varepsilon }(x)=\varepsilon x$ and $g_{\varepsilon
}=f_{\varepsilon }+n_{\varepsilon }$ where we have defined $n_{\varepsilon }$
as an increasing $\mathrm{C}^{\infty }$ function satisfying 
\begin{equation*}
n_{\varepsilon }(x)=%
\begin{cases}
-1 & \text{ if }x<-\frac{2}{\varepsilon }, \\ 
0 & \text{ if }x\in \left[ -\frac{1}{\varepsilon },\frac{1}{\varepsilon }%
\right] , \\ 
1 & \text{ if }x>\frac{2}{\varepsilon }.%
\end{cases}%
\end{equation*}%
Note that $\left( n_{\varepsilon }\right) _{\varepsilon }\in \mathcal{N}(%
\mathbb{R})$ as for any compact $K$, it exists $\varepsilon \in (0,1)$ such
that $n_{\varepsilon \mid K}\equiv 0$. But it can easily checked that $%
n_{\varepsilon }\notin \mathcal{N}_{\tau }(\mathbb{R})$. As $f_{\varepsilon
} $ is strictly increasing then $g_{\varepsilon }=f_{\varepsilon
}+n_{\varepsilon }$ is also. We have $f_{\varepsilon }^{-1}(y)=y/\varepsilon 
$ so $(f_{\varepsilon }^{-1})_{\varepsilon }\in \mathcal{X}_{\tau }(\mathbb{R%
})$, moreover it is easy to choose $n_{\varepsilon }$ so that $\left(
g_{\varepsilon }^{-1}\right) _{\varepsilon }\in \mathcal{X}_{\tau }(\mathbb{R%
})$. Now we have 
\begin{equation*}
g_{\varepsilon }^{-1}(y)=%
\begin{cases}
(y-1)/\varepsilon & \text{ if }y>3 \\ 
y/\varepsilon & \text{ if }y\in \lbrack -1,1].%
\end{cases}%
\end{equation*}%
So that $\left\vert \left\vert f_{\varepsilon }^{-1}-g_{\varepsilon
}^{-1}\right\vert \right\vert _{[3,4]}=1/\varepsilon $ which proves that $%
\left( f_{\varepsilon }^{-1}-g_{\varepsilon }^{-1}\right) _{\varepsilon
}\notin \mathcal{N}(\mathbb{R})$ .\newline
Moreover if we now turn back to the wave equation and we set $\psi (x)=x$
and $F=0$, the 2 solutions $(u_{\varepsilon })_{\varepsilon }$ and $%
(v_{\varepsilon })_{\varepsilon }$ corresponding to $(f_{\varepsilon
})_{\varepsilon }$ and $(g_{\varepsilon })_{\varepsilon }$ respectively are
given by 
\begin{equation*}
\begin{cases}
u_{\varepsilon }\left( x,y\right) =\varphi (x)+\chi _{\varepsilon
}^{f}(y)-\chi _{\varepsilon }^{f}(f_{\varepsilon }(x)), \\ 
v_{\varepsilon }\left( x,y\right) =\varphi (x)+\chi _{\varepsilon
}^{g}(y)-\chi _{\varepsilon }^{g}(g_{\varepsilon }(x)),%
\end{cases}%
\end{equation*}
\ so that we have 
\begin{equation*}
\frac{\partial (u_{\varepsilon }-v_{\varepsilon })}{\partial y}=\psi
(f_{\varepsilon }^{-1}(y))-\psi (g_{\varepsilon }^{-1}(y))=f_{\varepsilon
}^{-1}(y)-g_{\varepsilon }^{-1}(y).
\end{equation*}%
So $\left( \frac{\partial (u_{\varepsilon }-v_{\varepsilon })}{\partial y}%
\right) _{\varepsilon }\notin \mathcal{N}(\mathbb{R})$. This proves that the
usual equivalence is too coarse.
\end{remark}

Before proving the theorem we need the following
\begin{lemma}
\label{equiv_reciproc} Let $\left(f_{\e}\right)_\e,\left(g_{\e}\right)_\e \in \mathcal{X}_{\tau }(\mathbb{R})$ such that for every $\e $, $f_{\e},g_\e$ are bijective and $\left( f_{\e}^{-1}\right)_\e,\left( g_{\e}^{-1}\right)_{\e }\in \mathcal{X}_\tau(\mathbb{R})$. If moreover $\left( g_{\varepsilon }-f_{\varepsilon }\right)
_{\varepsilon }\in \mathcal{N}_{\tau }(\mathbb{R})$ we have that 
\begin{equation*}
\left(
f_{\varepsilon }^{-1}-g_{\varepsilon }^{-1}\right) _{\varepsilon }\in 
\mathcal{N}_{\tau }(\mathbb{R})
\end{equation*}
\end{lemma}
The proof will use the pointvalues characterization; so let us first define the following map (cf \cite{GKOS}, 1.2)
\begin{align*}
\Theta:\mathcal G_\tau(\R) &\to \mathcal F(\Rs)\\
\left[f_\e\right] &\mapsto f:\tilde x=[x_\e] \mapsto f(\tilde x)=[f_\e(x_\e)]
\end{align*}
where $\Rs$ denotes the field of generalized real numbers and $\mathcal F(\Rs)$ is the set of map from $\Rs$ to $\Rs$.
\begin{proof}
First $\mathcal G_\tau(\R)$ and $\mathcal F(\Rs)$ can be endowed with a structure of unitary rings where the operations are addition and composition of functions (the unit is then the identity function). Let us prove that $\Theta$ is a morphism between these rings; let $(f_\e)_\e,(g_\e)_\e \in \mathcal X_\tau(\R)$ and $f=\Theta([f_\e])$, $g=\Theta([g_\e])$, $h=\Theta([f_\e \circ g_\e])$ (note that $(f_\e \circ g_\e)_\e \in \mathcal X_\tau(\R)$):
\begin{align*}
\forall \tilde x=[x_\e] \in \Rs, h(\tilde x)=[(f_\e \circ g_\e)(x_\e)] &= [f_\e(g_\e(x_\e))]\\
&=f(g(\tilde x))\\
&=\Theta([f_\e]) \circ \Theta([g_\e]) (\tilde x)
\end{align*}
If we assume moreover that $f_{\e},g_\e$ are bijective, $\left( f_{\e}^{-1}\right),\left( g_{\e}^{-1}\right)_{\e }\in \mathcal{X}_\tau(\mathbb{R})$ and $\left( g_{\varepsilon }-f_{\varepsilon }\right)
_{\varepsilon }\in \mathcal{N}_{\tau }(\mathbb{R})$, we have that:
\begin{gather*}
\Id=[f_\e \circ f^{-1}_\e]= \Theta([f_\e]) \circ \Theta([f^{-1}_\e])\\
\Id=[f^{-1}_\e \circ f_\e]= \Theta([f^{-1}_\e]) \circ \Theta([f_\e])
\end{gather*}
So that $\Theta([f^{-1}_\e])=\Theta([f_\e])^{-1}$. Now as $\left[g_\e\right]=\left[f_\e\right]$, we have that $f=g$ so that $f^{-1}=g^{-1}$ and then $[f^{-1}_\e]=[g^{-1}_\e]$ which concludes the lemma.
\end{proof}

Now we can prove theorem \ref{thm_solution_unique}.

We have 
\begin{equation*}
u_{\varepsilon }(x,y)=u_{0,\varepsilon }(x,y)-{\iint\nolimits_{D(x,y,f_{%
\varepsilon })}}F(\xi ,\eta ,u_{\varepsilon }(\xi ,\eta ))d\xi d\eta ,
\end{equation*}%
where $u_{0,\varepsilon }(x,y)=\chi _{\varepsilon }^{f}\left( y\right) -\chi
_{\varepsilon }^{f}\left( f_{\varepsilon }\left( x\right) \right) +\varphi
(x)$ and $\chi _{\varepsilon }^{f}$ is a primitive of $\psi \circ
f_{\varepsilon }^{-1}$.\newline
So we take $\left( g_{\varepsilon }\right) _{\varepsilon }\in \mathcal{X}%
_{\tau }(\mathbb{R})$, such that $\left( f_{\varepsilon }-g_{\varepsilon
}\right) _{\varepsilon }\in \mathcal{N}_{\tau }(\mathbb{R})$; let $v=\left[
v_{\varepsilon }\right] $ the corresponding generalized solution. Let us
prove that $u=v$.\newline
We will in fact have prove a slightly stronger statement for reasons that
will be clear in the course of the proof%
\begin{equation*}
\forall K\Subset \mathbb{R}^{2},\forall \alpha \in \N^{2},\left(
P_{K_{\varepsilon },\alpha }(u_{\varepsilon }-v_{\varepsilon })\right)
_{\varepsilon }\in I_{A}
\end{equation*}%
Let us now fix $K\mathbb{\Subset }\mathbb{R}^{2}$. We have 
\begin{equation*}
v_{\varepsilon }(x,y)=v_{0,\varepsilon }(x,y)-{\iint\nolimits_{D(x,y,g_{%
\varepsilon })}}F(\xi ,\eta ,v_{\varepsilon }(\xi ,\eta ))d\xi d\eta \text{,}
\end{equation*}%
where $v_{0,\varepsilon }(x,y)=\chi _{\varepsilon }^{g}\left( y\right) -\chi
_{\varepsilon }^{g}(g_{\varepsilon }\left( x\right) )+\varphi (x)$ and $\chi
_{\varepsilon }^{g}$ is a primitive of $\psi \circ g_{\varepsilon }^{-1}$.%
\newline
So we get 
\begin{equation*}
u_{0,\varepsilon }(x,y)-v_{0,\varepsilon }(x,y)=\chi _{\varepsilon
}^{f}(y)-\chi _{\varepsilon }^{g}(y)-\chi _{\varepsilon }^{f}\left(
f_{\varepsilon }\left( x\right) \right) +\chi _{\varepsilon
}^{g}(g_{\varepsilon }\left( x\right) ).
\end{equation*}%
We compute 
\begin{equation*}
\dfrac{\partial }{\partial x}\left( -\chi _{\varepsilon }^{f}\left(
f_{\varepsilon }\left( x\right) \right) +\chi _{\varepsilon
}^{g}(g_{\varepsilon }\left( x\right) )\right) =\psi (x)\left(
f_{\varepsilon }^{\prime }(x)-g_{\varepsilon }^{\prime }(x)\right) \in 
\mathcal{N}_{\tau }(\mathbb{R}),
\end{equation*}%
but this implies that $\left\Vert \dfrac{\partial }{\partial x}\left( -\chi
_{\varepsilon }^{f}\left( f_{\varepsilon }\left( x\right) \right) +\chi
_{\varepsilon }^{g}(g_{\varepsilon }\left( x\right) )\right) \right\Vert
_{K_{\varepsilon }}\in I_{A}$. Indeed we can find $p\in \N$ such
that for any $m\in \N$ we have%
\begin{equation*}
\forall x\in \mathbb{R},\text{ }\dfrac{\partial }{\partial x}\left( -\chi
_{\varepsilon }^{f}\left( f_{\varepsilon }\left( x\right) \right) +\chi
_{\varepsilon }^{g}(g_{\varepsilon }\left( x\right) )\right) \leq
\varepsilon ^{m}(1+\left\vert x\right\vert )^{p},
\end{equation*}%
so 
\begin{equation*}
\left\Vert \dfrac{\partial }{\partial x}\left( -\chi _{\varepsilon
}^{f}\left( f_{\varepsilon }\left( x\right) \right) +\chi _{\varepsilon
}^{g}(g_{\varepsilon }\left( x\right) )\right) \right\Vert _{K_{\varepsilon
}}\leq \varepsilon ^{m}\max \left\{ (1+\left\vert a_{K,\varepsilon
}/2\right\vert )^{p}\right\} ,
\end{equation*}%
but $\left( (1+\left\vert a_{K,\varepsilon }/2\right\vert )^{p}\right)
_{\varepsilon }\in \left\vert A\right\vert $ thus we have obtained that 
\begin{equation*}
\left\Vert \dfrac{\partial }{\partial x}\left( -\chi _{\varepsilon
}^{f}\left( f_{\varepsilon }\left( x\right) \right) +\chi _{\varepsilon
}^{g}(g_{\varepsilon }\left( x\right) )\right) \right\Vert _{K_{\varepsilon
}}\in I_{A}.
\end{equation*}%
Then we obtain that $P_{K_{\varepsilon },1}\left( -\chi _{\varepsilon
}^{f}(f_{\varepsilon }(x))+\chi _{\varepsilon }^{g}(g_{\varepsilon
}(x))\right) \in I_{A}$. The proof is similar for higher derivatives.\newline
Now as $\left( f_{\varepsilon }^{-1}-g_{\varepsilon }^{-1}\right)
_{\varepsilon }\in \mathcal{N}(\mathbb{R})$ and $\psi \in \mathcal{O}_{M}$
then $\left( \psi \circ f_{\varepsilon }^{-1}-\psi \circ g_{\varepsilon
}^{-1}\right) _{\varepsilon }\in \mathcal{N}(\mathbb{R})$ and then we
finally obtain that 
\begin{equation*}
\forall \alpha ,\left( P_{K_{\varepsilon },\alpha }(u_{0,\varepsilon
}-v_{0,\varepsilon })\right) _{\varepsilon }\in I_{A}.
\end{equation*}%
We compute 
\begin{align*}
& u_{1,\varepsilon }(x,y)-v_{1,\varepsilon }(x,y) \\
=& \int_{f_{\varepsilon }^{-1}(y)}^{x}\int_{y}^{f_{\varepsilon }(\xi )}F(\xi
,\eta ,u_{\varepsilon }(\xi ,\eta ))d\eta d\xi -\int_{g_{\varepsilon
}^{-1}(y)}^{x}\int_{y}^{g_{\varepsilon }(\xi )}F(\xi ,\eta ,v_{\varepsilon
}(\xi ,\eta ))d\eta d\xi  \\
=& \int_{f_{\varepsilon }^{-1}(y)}^{x}\int_{y}^{f_{\varepsilon }(\xi )}\left[
F(\xi ,\eta ,u_{\varepsilon }(\xi ,\eta ))-F(\xi ,\eta ,v_{\varepsilon }(\xi
,\eta ))\right] d\eta d\xi  \\
& -\int_{g_{\varepsilon }^{-1}(y)}^{f_{\varepsilon
}^{-1}(y)}\int_{y}^{g_{\varepsilon }(\xi )}F(\xi ,\eta ,v_{\varepsilon }(\xi
,\eta ))d\eta d\xi -\int_{f_{\varepsilon }^{-1}(y)}^{x}\int_{f_{\varepsilon
}(\xi )}^{g_{\varepsilon }(\xi )}F(\xi ,\eta ,v_{\varepsilon }(\xi ,\eta
))d\eta d\xi .
\end{align*}%
As $f_{\varepsilon }\circ g_{\varepsilon }^{-1}\equiv id\mod \mathcal{N}%
^{s}(\mathbb{R})$, we have 
\begin{align*}
\sup_{y\in \lbrack -b,b]}\left\vert \int_{f_{\varepsilon
}^{-1}(y)}^{g_{\varepsilon }^{-1}(y)}\int_{y}^{g_{\varepsilon }(\xi )}F(\xi
,\eta ,v_{\varepsilon }(\xi ,\eta ))d\eta d\xi \right\vert & \leq
2b\int_{f_{\varepsilon }^{-1}(y)}^{g_{\varepsilon }^{-1}(y)}\sup_{\eta \in
\lbrack -b,b]}\left\vert F(\xi ,\eta ,v_{\varepsilon }(\xi ,\eta
))\right\vert d\xi  \\
& \leq 2b\left\Vert f_{\varepsilon }^{-1}-g_{\varepsilon }^{-1}\right\Vert
_{[-b,b]}\left\Vert F\right\Vert _{[\lambda _{\varepsilon },\mu
_{\varepsilon }]\times \lbrack -b,b]\times \mathbb{R}}
\end{align*}%
where 
\begin{equation*}
\begin{cases}
\lambda _{\varepsilon }=\min \{f_{\varepsilon }^{-1}(-b),g_{\varepsilon
}^{-1}(-b)\} \\ 
\mu _{\varepsilon }=\max \{f_{\varepsilon }^{-1}(b),g_{\varepsilon
}^{-1}(b)\},%
\end{cases}%
\end{equation*}%
which is negligible as $f_{\varepsilon }^{-1}\sim g_{\varepsilon }^{-1}$ and 
$\left\Vert F\right\Vert _{[\lambda _{\varepsilon },\mu _{\varepsilon
}]\times \lbrack -b,b]\times \mathbb{R}}\in \left\vert A\right\vert $.%
\newline
For the first derivative we have 
\begin{multline*}
\frac{d}{dy}\left( \int_{f_{\varepsilon }^{-1}(y)}^{g_{\varepsilon
}^{-1}(y)}\int_{y}^{g_{\varepsilon }(\xi )}F(\xi ,\eta ,v_{\varepsilon }(\xi
,\eta ))d\eta d\xi \right) =-\int_{y}^{g_{\varepsilon }(f_{\varepsilon
}^{-1}(y))}F(f_{\varepsilon }^{-1}(y),\eta ,v_{\varepsilon }(f_{\varepsilon
}^{-1}(y),\eta ))d\eta  \\
+\int_{f_{\varepsilon }^{-1}(y)}^{g_{\varepsilon }^{-1}(y)}F(\xi
,y,v_{\varepsilon }(\xi ,y))d\xi .
\end{multline*}%
And the same kind of arguments take care of those 2 terms. Now for the
higher derivatives 
\begin{align}
& \frac{d^{2}}{dy^{2}}\left( \int_{f_{\varepsilon }^{-1}(y)}^{g_{\varepsilon
}^{-1}(y)}\int_{y}^{g_{\varepsilon }(\xi )}F(\xi ,\eta ,v_{\varepsilon }(\xi
,\eta ))d\eta d\xi \right)   \label{AllD3} \\
=& -F(f_{\varepsilon }^{-1}(y),g_{\varepsilon }(f_{\varepsilon
}^{-1}(y)),v_{\varepsilon }(f_{\varepsilon }^{-1}(y),g_{\varepsilon
}(f_{\varepsilon }^{-1}(y))))+F(f_{\varepsilon }^{-1}(y),y,v_{\varepsilon
}(f_{\varepsilon }^{-1}(y),y))  \label{AllD4} \\
-& \int_{y}^{g_{\varepsilon }(f_{\varepsilon }^{-1}(y))}\frac{d}{dy}\left(
F(f_{\varepsilon }^{-1}(y),\eta ,v_{\varepsilon }(f_{\varepsilon
}^{-1}(y),\eta ))\right) d\eta   \label{AllD5} \\
+& F(g_{\varepsilon }^{-1}(y),y,v_{\varepsilon }(g_{\varepsilon
}^{-1}(y),y))-F(f_{\varepsilon }^{-1}(y),y,v_{\varepsilon }(f_{\varepsilon
}^{-1}(y),y))  \label{AllD6} \\
+& \int_{f_{\varepsilon }^{-1}(y)}^{g_{\varepsilon }^{-1}(y)}\frac{d}{dy}%
\left( F(\xi ,y,v_{\varepsilon }(\xi ,y))\right) d\xi .  \label{AllD7}
\end{align}%
The hypotheses on $F$ and the fact that $(g_{\varepsilon
}^{-1}-f_{\varepsilon }^{-1})_{\varepsilon }\in \mathcal{N}_{\tau }(\mathbb{R%
})$ takes care of the terms of lines (\ref{AllD4}) and (\ref{AllD6}). Let us
now turn our attention to line (\ref{AllD5}). We have 
\begin{multline*}
\int_{y}^{g_{\varepsilon }(f_{\varepsilon }^{-1}(y))}\frac{d}{dy}\left(
F(f_{\varepsilon }^{-1}(y),\eta ,v_{\varepsilon }(f_{\varepsilon
}^{-1}(y),\eta ))\right) d\eta  \\
=\dint_{y}^{g_{\varepsilon }(f_{\varepsilon }^{-1}(y))}[(f_{\varepsilon
}^{-1})^{\prime }(y)\frac{\partial F}{\partial \xi }(f_{\varepsilon
}^{-1}(y),\eta ,v_{\varepsilon }(f_{\varepsilon }^{-1}(y),\eta )) \\
+(f_{\varepsilon }^{-1})^{\prime }(y)\frac{\partial v_{\varepsilon }}{%
\partial x}(f_{\varepsilon }^{-1}(y),\eta )\frac{\partial F}{\partial z}%
(f_{\varepsilon }^{-1}(y),\eta ,v_{\varepsilon }(f_{\varepsilon
}^{-1}(y),\eta ))]d\eta .
\end{multline*}%
As $(g_{\varepsilon }\circ f_{\varepsilon }^{-1}-id)_{\varepsilon }\in 
\mathcal{N}_{\tau }(\mathbb{R})$ we can find a compact $L\subset \mathbb{R}$
such that 
\begin{equation*}
\forall \varepsilon ,\text{ }\{g_{\varepsilon }\circ f_{\varepsilon
}^{-1}(y):y\in \lbrack -b,b]\}\cup \lbrack -b,b]\subset L,
\end{equation*}%
and moreover it is sufficient to prove that 
\begin{multline*}
(\sup_{y\in \lbrack -b,b],\eta \in L}\mid (f_{\varepsilon }^{-1})^{\prime
}(y)\frac{\partial F}{\partial \xi }(f_{\varepsilon }^{-1}(y),\eta
,v_{\varepsilon }(f_{\varepsilon }^{-1}(y),\eta )) \\
+(f_{\varepsilon }^{-1})^{\prime }(y)\frac{\partial v_{\varepsilon }}{%
\partial x}(f_{\varepsilon }^{-1}(y),\eta )\frac{\partial F}{\partial z}%
(f_{\varepsilon }^{-1}(y),\eta ,v_{\varepsilon }(f_{\varepsilon
}^{-1}(y),\eta ))\mid )_{\varepsilon }\in \left\vert A\right\vert .
\end{multline*}%
But it is easy to see that 
\begin{equation*}
\left( \sup_{y\in \lbrack -b,b],\eta \in L}\left\vert (f_{\varepsilon
}^{-1})^{\prime }(y)\frac{\partial F}{\partial \xi }(f_{\varepsilon
}^{-1}(y),\eta ,v_{\varepsilon }(f_{\varepsilon }^{-1}(y),\eta ))\right\vert
\right) _{\varepsilon }\in \left\vert A\right\vert .
\end{equation*}%
For the other term the only part needing some new explanations is to prove
that 
\begin{equation*}
\left( \sup_{y\in \lbrack -b,b],\eta \in L}\frac{\partial v_{\varepsilon }}{%
\partial x}{(f_{\varepsilon }^{-1}(y),\eta )}\right) _{\varepsilon }\in
\left\vert A\right\vert .
\end{equation*}%
But here we use the fact that $(g_{\varepsilon }^{-1}-f_{\varepsilon
}^{-1})_{\varepsilon }\in \mathcal{N}_{\tau }(\mathbb{R})$ to find $%
\varepsilon _{0}$ such that 
\begin{equation}
\forall 0<\varepsilon <\varepsilon _{0},\left\Vert f_{\varepsilon
}^{-1}-g_{\varepsilon }^{-1}\right\Vert _{[-b,b]}<1.
\label{g_inv_f_inv_proches}
\end{equation}%
We proved in the proof of theorem \eqref{thm_solution_moderee} that 
\begin{equation*}
\left( P_{K_{\varepsilon },\alpha }(v_{\varepsilon })\right) _{\varepsilon
}\in \left\vert A\right\vert 
\end{equation*}%
and because of (\ref{g_inv_f_inv_proches}) we have that $f_{\varepsilon
}^{-1}(L)\times \lbrack -b,b]\subset K_{\varepsilon }$, which settles this
case.\newline
For higher derivatives the reasoning involves the same estimate and presents
no new obstacles. So this proves that 
\begin{equation*}
\forall \alpha ,\text{\thinspace }\left( P_{K_{\varepsilon },\alpha }\left(
\int_{f_{\varepsilon }^{-1}(y)}^{g_{\varepsilon
}^{-1}(y)}\int_{y}^{g_{\varepsilon }(\xi )}F(\xi ,\eta ,v_{\varepsilon }(\xi
,\eta ))d\eta d\xi \right) \right) _{\varepsilon }\in I_{A}.
\end{equation*}%
Similar arguments apply to prove that 
\begin{equation*}
\forall \alpha ,\,\left( P_{K_{\varepsilon },\alpha }\left(
\int_{f_{\varepsilon }^{-1}(y)}^{x}\int_{f_{\varepsilon }(\xi
)}^{g_{\varepsilon }(\xi )}F(\xi ,\eta ,v_{\varepsilon }(\xi ,\eta ))d\eta
d\xi \right) \right) _{\varepsilon }\in I_{A}.
\end{equation*}%
So we have proved that 
\begin{multline*}
P_{K_{\varepsilon },\alpha }\left( u_{1,\varepsilon }(x,y)-v_{1,\varepsilon
}(x,y)\right) _{\varepsilon } \\
\equiv \left( \int_{f_{\varepsilon }^{-1}(y)}^{x}\int_{y}^{f_{\varepsilon
}(\xi )}\left[ F(\xi ,\eta ,u_{\varepsilon }(\xi ,\eta ))-F(\xi ,\eta
,v_{\varepsilon }(\xi ,\eta ))\right] d\eta d\xi \right) _{\varepsilon }%
\mod I_{A}.
\end{multline*}%
We define 
\begin{equation*}
\sigma _{\varepsilon }(x,y)=u_{\varepsilon }(x,y)-v_{\varepsilon
}(x,y)-\int_{f_{\varepsilon }^{-1}(y)}^{x}\int_{y}^{f_{\varepsilon }(\xi )}
\left[ F(\xi ,\eta ,u_{\varepsilon }(\xi ,\eta ))-F(\xi ,\eta
,v_{\varepsilon }(\xi ,\eta ))\right] d\eta d\xi .
\end{equation*}%
So by the above arguments we just proved that $P_{K_{\varepsilon },\alpha
}\left( \sigma _{\varepsilon }\right) _{\varepsilon }\in I_{A}$. We now
define $w_{\varepsilon }(x,y)=u_{\varepsilon }(x,y)-v_{\varepsilon }(x,y)$.
Keeping the same notations as in the proof of theorem %
\eqref{thm_solution_moderee}, we want to prove that 
\begin{equation*}
\forall n,\,\left( P_{K_{\varepsilon },n}\left( w_{\varepsilon }\right)
\right) _{\varepsilon }\in I_{A}.
\end{equation*}%
Let us first prove that $P_{K_{\varepsilon },0}(w_{\varepsilon })\in I_{A}$.
First we have 
\begin{equation*}
F(\xi ,\eta ,u_{\varepsilon }(\xi ,\eta ))-F(\xi ,\eta ,v_{\varepsilon }(\xi
,\eta ))=w_{\varepsilon }(\xi ,\eta )\int_{0}^{1}\frac{\partial F}{\partial z%
}(\xi ,\eta ,u_{\varepsilon }(\xi ,\eta )+\theta (w_{\varepsilon }(\xi ,\eta
))d\theta ,
\end{equation*}%
then 
\begin{equation*}
w_{\varepsilon }(x,y)=\sigma _{\varepsilon }(x,y)+\int_{f_{\varepsilon
}^{-1}(y)}^{x}\int_{y}^{f_{\varepsilon }(\xi )}w_{\varepsilon }(\xi ,\eta
)\left( \int_{0}^{1}\frac{\partial F}{\partial z}(\xi ,\eta ,u_{\varepsilon
}(\xi ,\eta )+\theta (w_{\varepsilon }(\xi ,\eta ))d\theta \right) d\eta
d\xi .
\end{equation*}%
Now we have 
\begin{equation*}
\forall \varepsilon ,\,\cup _{(x,y)\in K_{\varepsilon }}\{(\xi ,\eta )\mid
\xi \in \lbrack x,f_{\varepsilon }^{-1}(y)],y\leq \eta \leq f_{\varepsilon
}(\xi )\}\subset L_{\varepsilon }=[\alpha _{K,\varepsilon },\beta
_{K,\varepsilon }]\times \lbrack -b,b],
\end{equation*}%
so that setting $l_{\varepsilon }=\sup_{L_{\varepsilon }\times \mathbb{R}%
}\left\vert \frac{\partial F}{\partial z}\right\vert $ we have 
\begin{equation*}
\left\vert w_{\varepsilon }(x,y)\right\vert \leq l_{\varepsilon
}\int_{\alpha _{K,\varepsilon }}^{\beta _{K,\varepsilon
}}\int_{y}^{f_{\varepsilon }(x)}\left\vert w_{\varepsilon }(\xi ,\eta
)\right\vert d\eta d\xi +\left\vert \sigma _{\varepsilon }(x,y)\right\vert ,
\end{equation*}%
so%
\begin{equation*}
\forall (x,y)\in K_{\varepsilon },\left\vert w_{\varepsilon
}(x,y)\right\vert \leq l_{\varepsilon }\int_{\alpha _{K,\varepsilon
}}^{\beta _{K,\varepsilon }}\int_{y}^{f_{\varepsilon }(x)}\left\vert
w_{\varepsilon }(\xi ,\eta )\right\vert d\eta d\xi +\left\Vert \sigma
_{\varepsilon }\right\Vert _{K_{\varepsilon }}.
\end{equation*}%
Letting $e_{\varepsilon }(y)=\sup_{\xi \in \lbrack \alpha _{K,\varepsilon
},\beta _{K,\varepsilon }]}\left\vert w_{\varepsilon }(\xi ,y)\right\vert $
we obtain 
\begin{equation*}
e_{\varepsilon }(y)\leq a_{K,\varepsilon }l_{\varepsilon
}\int_{y}^{b}e_{\varepsilon }(\eta )d\eta +\left\Vert \sigma _{\varepsilon
}\right\Vert _{K_{\varepsilon }}.
\end{equation*}%
So applying Gronwall's lemma we finally obtain 
\begin{equation*}
e_{\varepsilon }(y)\leq \left\Vert \sigma _{\varepsilon }\right\Vert
_{K_{\varepsilon }}\exp \left( a_{K,\varepsilon }l_{\varepsilon
}(b-y)\right) .
\end{equation*}%
Then 
\begin{equation*}
\forall (x,y)\in K_{\varepsilon },_{\,}\left\vert w_{\varepsilon
}(x,y)\right\vert \leq \left\Vert \sigma _{\varepsilon }\right\Vert
_{K_{\varepsilon }}\exp \left( a_{K,\varepsilon }l_{\varepsilon }2b\right) ..
\end{equation*}%
Consequently $P_{K_{\varepsilon },0}(w_{\varepsilon })\in I_{A}$. Which
implies the 0th order estimate. According to Proposition \ref{0-estimate},
we deduce $(w_{\varepsilon })_{\varepsilon }\in \mathcal{N}(\mathbb{R}^{2})$%
; consequently $u$ does not depend on the choice of the representative $%
\left( f_{\varepsilon }\right) _{\varepsilon }$ of the class $f=\left[
f_{\varepsilon }\right] \in \mathcal{G}_{\tau }\left( \mathbb{R}\right) $.

\bigskip

Now that the generalized solution is well defined let us prove that is is
indeed a generalization of the non-characteristic smooth case.\newline
From now on we will denote $\mathcal{X}^{s}$ and $\mathcal{N}^{s}$ the usual
special algebra and its ideal (that is the scale is polynomial in $%
\varepsilon $) and $\mathcal{G}^{s}=\mathcal{X}^{s}/\mathcal{N}^{s}$ the
usual special Colombeau algebra.\newline
Let us now see why theorems \ref{thm_solution_moderee} and \ref%
{thm_solution_unique} are a generalization of the classical case.\newline
Assume that the initial data are given along a path of equation $y=g(x)$
where $g$ is $\mathrm{C}^{\infty }$ and for all $x\in \mathbb{R}$, $%
g^{\prime }(x)\neq 0$, and that $\varphi ,\psi $ are $\mathrm{C}^{\infty }$.
We will then take $g_{\varepsilon }=g$ for all $\varepsilon $ to represent $g
$ in $\mathcal{G}^{s}$.\newline
Let $v$ be the classical solution then the generalized solution
corresponding to $g_{\varepsilon }$, is just $[v_{\varepsilon }]$ where $%
v_{\varepsilon }=v$ for all $\varepsilon $. Now what we need to prove is
that for any other representative $(f_{\varepsilon })_{\varepsilon }$ of $g$
in $\mathcal{A}^{s}$ we get a solution $(u_{\varepsilon })_{\varepsilon }$
which is moderate and equivalent to $(v_{\varepsilon })_{\varepsilon }$.
Note that this is not a consequence of theorems \ref{thm_solution_moderee}
and \ref{thm_solution_unique} as $g$ is not assumed to be tempered. But
actually we will see that the same proof works because $g^{-1}\in \mathrm{C}%
^{\infty }$ implies that $f_{\varepsilon }^{-1}$ is c-bounded. We will just
outline here the arguments.\newline
First because $f_{\varepsilon }^{-1}$ is c-bounded all bounds in (\ref{4})
are finite and then all sets there are also compact. Now looking at the
first arguments we do not need $\varphi ,\psi $ to be in $\mathcal{O}_{M}(%
\mathbb{R})$ anymore.

Moreover $\mathcal{C}=A/I_{A}$ is now overgenerated by $\left( \varepsilon
\right) _{\varepsilon }$, that is $\mathcal{A}=\mathcal{G}^{s}$. Now
stepping through the proof of Theorem \ref{thm_solution_moderee} we see that
because of the c-boundedness of $f_{\varepsilon }^{-1}$ the moderatness of $%
u_{0,\varepsilon }$ is obvious; recalling that $D_{\varepsilon }$\ is now
replaced by a fixed compact set, the following arguments go through easily,
proving the $0$th estimation for $\left( u_{\varepsilon }\right)
_{\varepsilon }$.

For the induction the c-boundedness of $f_{\varepsilon }^{-1}$ (and of
course of $\left( f_{\varepsilon }\right) _{\varepsilon }$) ensures that all
integral are done over bounded sets independent of $\varepsilon $ \ which
removes the need for temperatness to ensure that these integrals are
moderate.

For the independance proof the lemma \ref{equiv_reciproc}\ is not needed here
as both $f_{\varepsilon }$ and $f_{\varepsilon }^{-1}$ are c-bounded which
implies that composition works as expected without the functions being
tempered. The same remarks apply to the proof of Theorem \ref%
{thm_solution_unique} showing that it remains true in this context.

To summarize we have proved the

\begin{theorem}
The solution to the Cauchy problem ($P_{C})$ in the classical case when the
initial data are smooth and given along the curve $y=f(x)$ with $f\in 
\mathrm{C}^{\infty }(\mathbb{R})$ and $f^{\prime }(x)\neq 0$ for all $x\in 
\mathbb{R}$, coincides with any generalized solution associated to any $%
(f_{\varepsilon })_{\varepsilon }\in \mathcal{A}^{s}(\mathbb{R})$ such that $%
[f]=[f_{\varepsilon }]\in \mathcal{A}^{s}(\mathbb{R}).$
\end{theorem}

\section{Examples}

\label{sect_exemples} We now compute 2 examples where our method provides us
with a generalized solution in characteristic situations.

\begin{example}
We consider the characteristic Cauchy problem where the initial
values are smooth functions given on the characteristic curve $\gamma $
whose equation is $y=x^{3}$. We approach $\gamma $
by $\gamma _{\varepsilon }$ whose equation is
$y=x^{3}+\varepsilon x$. We suppose that $F=0.$ According to
the previous notations, we have to solve
\begin{equation*}
P_{1,\varepsilon}
\begin{cases}
\dfrac{\partial ^{2}u_{\varepsilon }}{\partial x\partial y}(x,y)=0 &(1)\\
u_{\varepsilon }\left( x,x^{3}+\varepsilon x\right) =\varphi (x) &(2)\\
\dfrac{\partial u_{\varepsilon }}{\partial y}\left( x,x^{3}+\varepsilon
x\right) =\psi (x) &(3)
\end{cases}
\end{equation*}
We can solve $P_{1,\varepsilon }$ by putting $u_{\varepsilon
}(x,y)=h_{\varepsilon }\left( x\right) +k_{\varepsilon }\left( y\right) $
and we have 
\begin{equation}
(2)\iff h_{\varepsilon }\left( x\right) +k_{\varepsilon }(x^{3}+\varepsilon
x)=\varphi (x),\text{ \ }(3)\iff k_{\varepsilon }^{\prime }\left(
x^{3}+\varepsilon x\right) =\psi (x)  \tag{V2}  \label{V2}
\end{equation}%
that is to say\textit{\ }$k_{\varepsilon }^{\prime }\left( t\right) =\psi
\left( f_{\varepsilon }^{-1}\left( t\right) \right) $ where $%
x^{3}+\varepsilon x=f_{\varepsilon }(x)=t$ is equivalent to$\
x=f_{\varepsilon }^{-1}\left( t\right) $. As $\varphi (x)=h_{\varepsilon
}\left( x\right) +k_{\varepsilon }(x^{3}+\varepsilon x)$, we deduce $\varphi
^{\prime }(x)=h_{\varepsilon }^{\prime }\left( x\right) +(3x^{2}+\varepsilon
)k_{\varepsilon }^{\prime }(x^{3}+\varepsilon x)$ then%
\begin{eqnarray*}
\tint\nolimits_{x}^{f_{\varepsilon }^{-1}\left( y\right) }\varphi ^{\prime
}(\xi )d\xi &=&\tint\nolimits_{x}^{f_{\varepsilon }^{-1}\left( y\right)
}h_{\varepsilon }^{\prime }(\xi )d\xi +\tint\nolimits_{x}^{f_{\varepsilon
}^{-1}\left( y\right) }\left( 3\xi ^{2}+\varepsilon \right) k_{\varepsilon
}^{\prime }(\xi ^{3}+\varepsilon \xi )d\xi \\
&=&\tint\nolimits_{x}^{f_{\varepsilon }^{-1}\left( y\right) }h_{\varepsilon
}^{\prime }(\xi )d\xi +\tint\nolimits_{x}^{f_{\varepsilon }^{-1}\left(
y\right) }\left( 3\xi ^{2}+\varepsilon \right) \psi (\xi )d\xi \\
&=&\left[ h_{\varepsilon }(\xi )\right] _{x}^{f_{\varepsilon }^{-1}\left(
y\right) }+\tint\nolimits_{0}^{f_{\varepsilon }^{-1}\left( y\right) }\left(
3\xi ^{2}+\varepsilon \right) \psi (\xi )d\xi -\tint\nolimits_{0}^{x}\left(
3\xi ^{2}+\varepsilon \right) \psi (\xi )d\xi .
\end{eqnarray*}%
Then%
\begin{equation*}
\varphi (f_{\varepsilon }^{-1}\left( y\right) )-\varphi \left( x\right)
=h_{\varepsilon }\left( f_{\varepsilon }^{-1}\left( y\right) \right)
-h_{\varepsilon }\left( x\right) +\tint\nolimits_{0}^{f_{\varepsilon
}^{-1}\left( y\right) }\left( 3\xi ^{2}+\varepsilon \right) \psi (\xi )d\xi
-\tint\nolimits_{0}^{x}\left( 3\xi ^{2}+\varepsilon \right) \psi (\xi )d\xi .
\end{equation*}%
We deduce that%
\begin{equation*}
\varphi (f_{\varepsilon }^{-1}\left( y\right) )-h_{\varepsilon }\left(
f_{\varepsilon }^{-1}\left( y\right) \right) +h_{\varepsilon }\left(
x\right) =\varphi \left( x\right) +\tint\nolimits_{0}^{f_{\varepsilon
}^{-1}\left( y\right) }\left( 3\xi ^{2}+\varepsilon \right) \psi (\xi )d\xi
-\tint\nolimits_{0}^{x}\left( 3\xi ^{2}+\varepsilon \right) \psi (\xi )d\xi .
\end{equation*}%
According to (\ref{V2}), we obtain $\varphi (f_{\varepsilon }^{-1}\left(
y\right) )-h_{\varepsilon }\left( f_{\varepsilon }^{-1}\left( y\right)
\right) =k_{\varepsilon }\left( f_{\varepsilon }\left( f_{\varepsilon
}^{-1}\left( y\right) \right) \right) =k_{\varepsilon }\left( y\right) $, so,%
\begin{equation*}
k_{\varepsilon }\left( y\right) +h_{\varepsilon }\left( x\right) =\varphi
\left( x\right) +\tint\nolimits_{0}^{f_{\varepsilon }^{-1}\left( y\right)
}\left( 3\xi ^{2}+\varepsilon \right) \psi (\xi )d\xi
-\tint\nolimits_{0}^{x}\left( 3\xi ^{2}+\varepsilon \right) \psi (\xi )d\xi ,
\end{equation*}%
that is to say,%
\begin{equation*}
u_{\varepsilon }(x,y)=\varphi (x)+\tint\nolimits_{0}^{f_{\varepsilon
}^{-1}\left( y\right) }\left( 3\xi ^{2}+\varepsilon \right) \psi (\xi )d\xi
-\tint\nolimits_{0}^{x}\left( 3\xi ^{2}+\varepsilon \right) \psi (\xi )d\xi .
\end{equation*}%
Solve $P_{1,\varepsilon }$. Set 
\begin{equation*}
H\left( x\right) =\varphi (x)-\tint\nolimits_{0}^{x}3\xi ^{2}\psi (\xi )d\xi
;\text{ }L\left( x\right) =\tint\nolimits_{0}^{x}\psi (\xi )d\xi .
\end{equation*}%
We can write 
\begin{equation*}
h_{\varepsilon }\left( x\right) =\varphi (x)-\tint\nolimits_{0}^{x}3\xi
^{2}\psi (\xi )d\xi +\varepsilon \tint\nolimits_{0}^{x}\psi (\xi )d\xi
=H\left( x\right) +\varepsilon L\left( x\right)
\end{equation*}%
and%
\begin{equation*}
k_{\varepsilon }\left( y\right) =\tint\nolimits_{0}^{f_{\varepsilon
}^{-1}\left( y\right) }\left( 3\xi ^{2}+\varepsilon \right) \psi (\xi )d\xi .
\end{equation*}%
So%
\begin{equation*}
u_{\varepsilon }=h_{\varepsilon }\otimes 1_{y}+1_{x}\otimes k_{\varepsilon }
\end{equation*}%
then 
\begin{equation*}
h_{\varepsilon }\otimes 1_{y}=H\otimes 1_{y}+\varepsilon L\otimes 1_{y}.
\end{equation*}%
The class of $\left( h_{\varepsilon }\otimes 1_{y}\right) _{\varepsilon }$
lies in a Colombeau algebra but it is not the case for $\left( 1_{x}\otimes
k_{\varepsilon }\right) _{\varepsilon }$ and we have to involve a convenient 
$(\mathcal{C},\mathcal{E},\mathcal{P})$-algebra.
\end{example}

\begin{example}
Let us take $f(x)=sgn(x)$, and 
\begin{equation*}
f_{\varepsilon }(x)=\tanh (x/\varepsilon )+\alpha _{\varepsilon }(x)\text{
and }g_{\varepsilon }(x)=\frac{2}{\pi }\arctan (x/\varepsilon )+\alpha
_{\varepsilon }(x),
\end{equation*}%
where $(\alpha _{\varepsilon })_{\varepsilon }$ is any family of smooth
functions with $\alpha _{\varepsilon }(\mathbb{R})=\mathbb{R}$ and $\alpha
_{\varepsilon }^{\prime }>0$ such that $(\alpha _{\varepsilon
})_{\varepsilon }\in \mathcal{N}$. One can check that $f_{\varepsilon
}\rightarrow f$ and $g_{\varepsilon }\rightarrow f$ (simple limits, and in
fact uniform convergence over compacts not containing $0$). Moreover $%
f_{\varepsilon },g_{\varepsilon }\in \mathcal{X}_{\tau }$ but $%
f_{\varepsilon }-g_{\varepsilon }\notin \mathcal{N}$. We compute $%
u_{\varepsilon }=u_{0,\varepsilon }$ as $F=0$, using both regularizations:
first for $f_{\varepsilon }$ we have $u_{\varepsilon }(x,y)=\chi
_{\varepsilon }(y)-\chi _{\varepsilon }(f_{\varepsilon }(x))+\varphi (x)$, so%
\begin{align*}
\frac{\partial {u_{\varepsilon }}}{\partial {x}}& =-f_{\varepsilon }^{\prime
}(x)\psi (x)+\varphi ^{\prime }(x) \\
& \equiv -\frac{1}{\varepsilon \cosh ^{2}(x/\varepsilon )}\psi (x)+\varphi
^{\prime }(x)\text{ mod }\mathcal{N},
\end{align*}%
then%
\begin{equation*}
\left( \frac{\partial {u_{\varepsilon }}}{\partial {x}}\right) _{\varepsilon
}\thicksim -2\psi (0)\delta +\varphi ^{\prime }(x)\text{ as }\int_{R}\tanh
(x)dx=2.
\end{equation*}%
We also have 
\begin{equation*}
\frac{\partial {u_{\varepsilon }}}{\partial {y}}=\psi (f_{\varepsilon
}^{-1}(y))\equiv \psi (\varepsilon \tanh ^{-1}y)\text{ mod }\mathcal{N},
\end{equation*}%
so%
\begin{equation*}
\left( \frac{\partial {u_{\varepsilon }}}{\partial {y}}\right) _{\varepsilon
}\underset{}{\underset{\mathrm{C}^{\infty }}{\longrightarrow }}\psi (0)
\end{equation*}%
Then we can conclude that 
\begin{equation*}
(u_{\varepsilon })_{\varepsilon }\thicksim -2\psi (0)Y_{x}+\varphi (x)+\psi
(0)y
\end{equation*}%
Similar computations give that $(v_{\varepsilon })_{\varepsilon }$ is
associated to the same distribution. As a side note this example shows that
nonequivalent (but associated) deformations of the characteristic curve give
generalized solutions associated to the same distribution. We do not know
yet if this is a more general phenomenon.
\end{example}

\end{document}